\documentclass[12pt]{amsart}
\usepackage{enumerate, amssymb}
\usepackage[colorlinks]{hyperref}
\input{diagrams.tex}
\def\pdfsyncstart{}
\def\pdfsyncstop{}

\def\bdi{\pdfsyncstop\begin{diagram}}
\def\edi{\end{diagram}\pdfsyncstart}

\theoremstyle{plain}

\newtheorem{thm}{Theorem}[section]
\newtheorem{cor}[thm]{Corollary}
\newtheorem{lem}[thm]{Lemma}
\newtheorem{prop}[thm]{Proposition}

\theoremstyle{definition}
\newtheorem{defi}[thm]{Definition}
\newtheorem{defis}[thm]{Definitions}
\newtheorem{conj}[thm]{Conjecture}
\newtheorem{conv}[thm]{Convention}
\newtheorem{nota}[thm]{Notation}
\newtheorem{rem}[thm]{Remark}
\newtheorem{rems}[thm]{Remarks}
\newtheorem{exa}[thm]{Example}
\newtheorem{exas}[thm]{Examples}
\newtheorem{sit}[thm]{}

\newcommand{\brem}{\begin{rem}}
\newcommand{\brems}{\begin{rems}}
\newcommand{\erem}{\end{rem}}
\newcommand{\erems}{\end{rems}}
\newcommand{\bexa}{\begin{exa}}
\newcommand{\bexas}{\begin{exas}}
\newcommand{\eexa}{\end{exa}}
\newcommand{\eexas}{\end{exas}}
\newcommand{\bdefi}{\begin{defi}}
\newcommand{\edefi}{\end{defi}}
\newcommand{\bdefis}{\begin{defis}}
\newcommand{\edefis}{\end{defis}}
\newcommand{\bcor}{\begin{cor}}
\newcommand{\ecor}{\end{cor}}
\newcommand{\blem}{\begin{lem}}
\newcommand{\elem}{\end{lem}}
\newcommand{\bconv}{\begin{conv}}
\newcommand{\econv}{\end{conv}}
\newcommand{\bconj}{\begin{conj}}
\newcommand{\econj}{\end{conj}}
\newcommand{\bprop}{\begin{prop}}
\newcommand{\eprop}{\end{prop}}
\newcommand{\bthm}{\begin{thm}}
\newcommand{\ethm}{\end{thm}}
\newcommand{\bnota}{\begin{nota}}
\newcommand{\enota}{\end{nota}}
\newcommand{\bsit}{\begin{sit}}
\newcommand{\esit}{\end{sit}}
\newcommand{\be}{\begin{eqnarray}}
\newcommand{\ee}{\end{eqnarray}}
\newcommand{\bproof}{\begin{proof}}
\newcommand{\eproof}{\end{proof}}
\def\ba{\begin{array}}
\def\ea{\end{array}}

\def\bnum{\begin{enumerate}}
\def\enum{\end{enumerate}}

\newcommand{\R}{{\mathbb R}}
\newcommand{\C}{{\mathbb C}}
\newcommand{\Q}{{\mathbb Q}}
\newcommand{\Z}{{\mathbb Z}}
\newcommand{\N}{{\mathbb N}}

\newcommand{\G}{{\Gamma}}

\newcommand{\p}{{\partial}}

\newcommand{\nlin}{\unitlength1mm\begin{picture}(0,9.25)
                       \put(0,0.75){\line(0,1){8.5}}
                      \end{picture}}
\newcommand{\slin}{\unitlength1mm\begin{picture}(0,0)
                       \put(0,-0.75){\line(0,-1){8.5}}
                      \end{picture}}
\newcommand{\nolin}{\unitlength1mm\begin{picture}(0,7)
                       \put(0.53,0.53){\line(1,1){8.94}}
                      \end{picture}}
\newcommand{\solin}{\unitlength1mm\begin{picture}(0,0)
                       \put(0.53,-0.53){\line(1,-1){8.94}}
                      \end{picture}}
\newcommand{\nwlin}{\unitlength1mm\begin{picture}(0,7)
                       \put(-0.53,0.53){\line(-1,1){8.94}}
                      \end{picture}}
\newcommand{\swlin}{\unitlength1mm\begin{picture}(0,0)
                       \put(-0.53,-0.53){\line(-1,-1){8.94}}
                      \end{picture}}
\newcommand{\vlin}[1]{\hspace{0.75mm}\unitlength1mm\begin{picture}(#1,0)
                       \put(0,0){\line(1,0){#1}}
                      \end{picture}\hspace{0.75mm}\rule[-3mm]{0mm}{4mm}}

\def\llin{\vlin{11.5}}
\newcommand{\lin}{\vlin{8.5}}

\newcommand{\co}[1]{\unitlength1mm\begin{picture}(0,8)
    \put(0,0){\circle{1.5}}
    \put(0,3){\makebox(0,5)[b]{$#1$}}
                      \end{picture}}

\newcommand{\cou}[2]{\unitlength1mm\begin{picture}(0,8)
    \put(0,0){\circle{1.5}}
    \put(0,3){\makebox(0,5)[b]{$#1$}}
    \put(0,-7){\makebox(0,4)[t]{$#2$}}
      \end{picture}
      \rule[-7mm]{0mm}{7mm}}

\newcommand{\cshiftup}[2]{\unitlength1mm\begin{picture}(0,9.25)
                       \put(0,10){\cou{#1}{#2}}
                      \end{picture}}

\newcommand{\cshiftdown}[2]{\unitlength1mm\begin{picture}(0,9.25)
                       \put(0,-10){\cou{#1}{#2}}
                      \end{picture}}

\title{Birational transformations of weighted graphs}

\author{Hubert Flenner}
\address{Fakult\"at f\"ur Mathematik,
Ruhr Universit\"at Bochum, Geb.\ NA 2/72, Universit\"ats\-str.\
150, 44780 Bochum, Germany}
\email{Hubert.Flenner@rub.de}

\author{Shulim Kaliman}
\address{Department of Mathematics,
University of Miami, Coral Gables, FL  33124, U.S.A.}
\email{kaliman@math.miami.edu}

\author{Mikhail Zaidenberg}
\address{Universit\'e
Grenoble I, Institut Fourier, UMR 5582 CNRS-UJF, BP 74, 38402
St.\ Martin d'H\`eres c\'edex, France}
\email{zaidenbe@ujf-grenoble.fr}

\thanks{
{\bf Acknowledgements:} This research was started during a visit
of the first two authors at the Fourier Institute, Grenoble, and
continued during a visit of two of us at the Max Planck Institute
of Mathematics in Bonn. We thank these institutions for their
generous support and excellent working conditions. The research
was also partially supported by NSA Grant no. MDA904-00-1-0016 and
the DFG-Schwerpunkt Komplexe Geometrie. }

\thanks{
\mbox{\hspace{11pt}}{\it 1991 Mathematics Subject Classification}:
05C99, 14R99, 14J50.\\
\mbox{\hspace{11pt}}{\it Key words}: Weighted graph, standard
model, minimal model, birational transformation, elementary
transformation}

\dedicatory{Dedicated to Masayoshi Miyanishi}
\begin{document}

\begin{abstract}
We introduce the notion of a standard weighted graph
and show that every weighted graph has an essentially unique standard model.
Moreover we classify birational transformations between such
models. Our central result shows that these are composed of
elementary transformations. The latter ones are defined similarly
to the well known elementary transformations of ruled surfaces.

In a forthcoming paper, we apply these results in the geometric
setup to obtain standard equivariant completions of affine
surfaces with an action of certain algebraic groups. We show that
these completions are unique up to equivariant elementary
transformations.
\end{abstract}

\maketitle

{\footnotesize \tableofcontents}

\section{Introduction}

Birational transformations of weighted graphs were studied by many authors,
mainly due to their importance for understanding completions
of algebraic surfaces, see e.g.\
\cite{DaGi, Dai1, Dai2, FlZa0, Fu, Hi2, Mu, Ne, Ra, Ru}.
Danilov and Gizatullin \cite{DaGi} were the first to introduce several
special forms of linear graphs like semistandard, $m$-standard or
quasistandard ones to deduce their interesting results on automorphism
groups of affine algebraic surfaces that  are completed by a chain of
rational curves.  More recently Daigle \cite{Dai1, Dai2} studied
standard models of weighted trees and showed that any such tree has
a unique standard model in its birational equivalence class.

In this paper we generalize this theory to arbitrary weighted graphs.
The Reduction Theorem
\ref{graph.7}
shows that any graph, possibly with cycles, loops
and multiple edges, admits a standard model.
Moreover, this standard model is essentially unique
in its birational equivalence class, see
Corollary \ref{unique}.

A major part of the paper is devoted to the study of birational
transformations of standard weighted graphs. Any such
transformation preserves the branching points, see Lemma
\ref{branching}. Therefore it is sufficient to classify them for
linear chains and for circular graphs. In both cases we show that
one can decompose such a transformation into simpler ones called
{\it moves, shifts and turns}, see Propositions \ref{p1}, \ref{p2}
and Theorem \ref{maincirc}. In particular (see  Theorem
\ref{graph.8}) any birational transformation between standard
models is composed of elementary transformations, which are
defined similar to those for ruled surfaces (see Definition
\ref{graph.5}).

It is worthwhile to compare our result with those in the paper of
Danilov and Gizatullin \cite{DaGi}. Indeed, Theorem 1 in
\cite{DaGi} implies in particular that any biregular map between
two affine surfaces completed by standard linear chains\footnote{I
Notice that in \cite{DaGi} a graph is called standard, if it is
semistandard, $m$-standard or quasistandard, whereas our standard
graphs are $0$-standard in the terminology of \cite{DaGi}.} of
rational curves, is a product of birational elementary
transformations and some standard reconstructions of good
completions of our surfaces. While proving this theorem, the
control on the indeterminacy points of the underlying birational
map was important. Since we deal with general graphs which include
non-simply connected ones and do not necessarily satisfy the Hodge
index theorem, we do not have an underlying birational map of
algebraic surfaces. Therefore, our pivotal point is different.
Given a birational transformation of two standard weighted graphs
$A$ and $B$, we look for dominant maps from a third graph $\G\to
A$ and $\G\to B$. Actually we decompose our birational
transformation into a sequence of elementary transformations
dominated at every step by $\G$ so that one can apply our results
to the corresponding situation of boundary divisors of algebraic
surfaces. The role of indeterminacy points in \cite{DaGi} is
played by the vertices of $\G$ that are not contracted in both $A$
and $B$.

To complete the picture, we survey in the Appendix some well known
facts on the adjacency matrix and the discriminant of a weighted
graph and their  behaviour under birational transformations. In
particular, we compute the spectra of standard weighted graphs.

In the subsequent paper \cite{FKZ}  we will apply our results in
the geometric setting to obtain equivariant standard completions
of affine surfaces equipped with an effective action of certain
algebraic groups, cf. \cite[\S 6]{DaGi}.

\section{
Weighted graphs}

\subsection{Generalities}
  A {\em (combinatorial) graph } consists of a nonempty set of
vertices $\G^{(0)}$ and a set of edges $\G^{(1)}$ together with a
boundary map $\p$ which associates to every edge
$e\in\G^{(1)}$ the set $\p(e)$
consisting of one or two vertices, called the {\em end points}
of $e$. An edge $e$ with
just one end point is a simple loop.
In this subsection we consider weighted graphs with arbitrary real
weights of vertices, and we denote by $|\G|$ the number of
vertices of $\G$. All our graphs are assumed to be finite.

The {\it degree} (or the {\it valency}) $\deg_\Gamma
(v)$ of a vertex $v\in \G^{(0)}$ is the number of edges adjacent
to $v$, where we count the loops at $v$ twice.  The {\it
branches} of a connected graph $\Gamma$ at $v$ are the
connected components of the
graph $\Gamma\ominus v$ obtained from $\Gamma$ by deleting the
vertex $v$ and all its incident edges. In case $\deg (v)>2$ we
call $v$ a {\em branching point}; if $\deg (v)\le 1$, $v$ is
called an {\em end vertex} or a {\it tip}, and a {\em linear vertex} if $\deg
(v)=2$.

A graph is said to be {\em linear} or a {\em chain} if it has two
end vertices and all other vertices are linear. By a {\em
circular} graph we mean a connected graph with only linear
vertices. We let $B(\G)$ denote the set of all branching points of
$\G$. A connected graph $\G$ with $B(\G)=\emptyset$ is either
linear or circular.

The connected components of $\G\ominus B(\G)$ will be called the
{\em segments} of $\G$. Clearly the segments of $\G$ are either
linear or circular weighted graphs, as they do not include the
branching points of $\G$. Moreover for a connected graph $\G$ a
circular segment can appear if and only if $\G$ is circular
itself.

The {\it branching number} at $v$ is $\nu_\Gamma (v)=\max \{0,
\deg (v)-2\}$, and the {\it total branching number} is
$$\nu(\Gamma)=\sum_{v\in B(\G)} \nu_\Gamma (v)\,.$$
\bsit\label{birtr} We denote by $[[w_0,\ldots, w_n]]$ a chain
with linearly ordered vertices $v_0<v_1<\ldots <v_n$, where
$w_i=v_i^2\in\R$ is the weight of $v_i$:
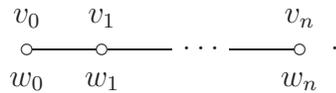
$$
\cou{v_0}{w_0}\lin \cou{v_1}{w_1}\lin \ldots
\lin\cou{v_n}{w_n}\quad.$$ Similarly, we denote by $((w_0,\ldots,
w_n))$ a circular graph with cyclically ordered vertices
$v_0<\ldots<v_n<v_0\ldots$ with weights $w_i=v_i^2\in\R$ of $v_i$.

\bsit\label{notation} Given two ordered linear chains
$L$ with vertices $v_1<\ldots<v_k$ and $M$  with vertices $u_1<\ldots<u_l$
we denote by $LM$ their {\it join}
that is, the ordered linear chain
 with vertices $v_1<\ldots<v_k<u_1<\ldots<u_l$. We let $L^{-1}$
be the chain $L$ with the reversed ordering $v_k<\ldots<v_1$.
For a sequence of ordered linear chains $L_1,\ldots,L_n$, we let
$((L_1\ldots L_n))$
be the circular cyclically ordered graph made of
their join.
\esit

\bdefis\label{graph.1} For a weighted graph $\Gamma$,
an {\it inner blowup} $\Gamma'\to \Gamma$ at an edge $e$ with end
vertices $v_0,v_1\in \G^{(0)}$ consists in
introducing a new vertex $v\in (\G')^{(0)}$ of weight $-1$
subdividing $e$ in two edges $e'$ and $e''$ with
$\p(e')=\{v_0,v\}$ and $\p(e'')=\{v,v_1\}$, and diminishing by
$1$ the weights of $v_0$ and $v_1$ (in case where $e$ is a loop
i.e., $v_0=v_1$, the weight of $v_0$ is diminishing by 2). An
{\it outer blowup} $\Gamma'\to \Gamma$ at a vertex $v_0$ of
$\Gamma$ consists in introducing a new vertex $v$ of weight $-1$
and a new edge $e$ with end vertices $v_0,v$, and diminishing by $1$
the weight of $v_0$. In both cases, the inverse procedure is
called {\it blowdown} of $v$.

The graph $\Gamma$ is {\it minimal} if it does not admit any
blowdown. Clearly $\Gamma$ is minimal if and only if every
segment of $\Gamma$ is.

A {\em birational transformation} of a graph $\Gamma$ into another
one $\Gamma'$ is a sequence of blowing ups and downs. We write
$\G\sim \G'$ or $\G\dasharrow \G'$ if such a transformation does
exist, and $\Gamma\to \Gamma'$ if $\Gamma$ is obtained from
$\Gamma'$ by a sequence of only blowups. In the latter case we say
that $\G$ {\it dominates} $\G'$, and we call $\Gamma\to \Gamma'$ a
{\it domination}. If $\Gamma\to \Gamma'$ is a domination and $v$
is a vertex of $\G'$, we denote by $\hat v$ the corresponding
vertex of $\G$ called the {\it proper transform} or the {\it
preimage} of $v$. Similarly, for a subgraph $A$ of $\G'$ with
vertices $\{a_i\}$ , $\hat A$ stands for a subgraph of $\G$ with
vertices $\{\hat a_i\}$.

Any birational transformation $\gamma:\Gamma_1\to \Gamma_2$
fits into a commutative diagram
\bdi[size=1.5em]
   &&&&& \Gamma \\
   &&&&\ldTo &&\rdTo \\
   &&& \Gamma_1&&\rDashto^{\gamma}&&\Gamma_2\,,
\edi where $\Gamma\to\Gamma_i$, $i=1,2$, are dominations.
Moreover we may suppose that this decomposition is {\it
relatively minimal} that is, no $(-1)$-vertex of $\Gamma$ is
contracted in both directions, see \cite[Appendix, Remark
A.1(1)
]{FlZa0}. \edefis

Clearly, the topological (homotopy) type
of a graph is birationally invariant whereas $\nu(\G)$ is not, in
general.

We recall the following facts, see \cite[Appendix to \S 4]{FlZa0}
and also \cite[Cor. 3.6]{Dai1}, \cite{Ru}.
We provide a short argument in our more general setting.

\blem\label{branching} Let $\Gamma$ dominate a minimal weighted
graph $\Gamma_1$.
\begin{enumerate}
\item[(a)] If a branching point $v$ of degree $r$ in $\Gamma$ is
not contained in $B(\Gamma_1)\hat{}$, then $\Gamma$ has at
least $r-2$ branches at $v$ contractible inside $\Gamma$ and
contracted in $\Gamma_1$.

\item[(b)] If $\Gamma\to\G_1$, $\G\to \Gamma_2$ is a pair of
dominations and $\G_1,\,\G_2$ are minimal then
$B(\Gamma_1)\hat{}=B(\Gamma_2)\hat{}$.
\end{enumerate} \elem

\bproof  To show (a) we note that a vertex $v\in B(\G)$ of degree
$r\ge 3$ can become at most linear only after contracting $r-2$
branches of $\G$ at $v$. Moreover, blowdowns in one of them do not
affect the other ones, so each of them must be contractible.

To show (b), suppose on the contrary that there exists a vertex
$\hat v_1\in B(\G_1)\hat{}\ominus B(\Gamma_2)\hat{}$, and let
$r\ge 3$ be the degree of $\hat v_1$ in $\G$. According to (a)
there are $r-2$ contractible branches of $\G$ at $\hat v_1$ that
are blown down in $\G_2$. Since $v_1\in B(\G_1)$, at least one of
these branches is not blown down completely in $\G_1$. This
contradicts the minimality assumption. \eproof

\bcor\label{brpts} The number of branching points of a minimal
graph $\G$, their degrees and the total branching number $\nu
(\G)$ are birational invariants. In particular, a weighted graph
$\Gamma$ that can be transformed into one with fewer branching
points is not minimal.\ecor


Thus the only birationally non-rigid elements of a minimal graph
can be its segments. A graph with no segments is birationally
rigid that is, has a birationally unique minimal model. A segment
of a graph can eventually be non-rigid even being minimal, see an
example in \ref{graph.5} below. However the number of such
segments and their types (linear or circular) remain stable under
birational transformations.

\bexas\label{nontip} The graph

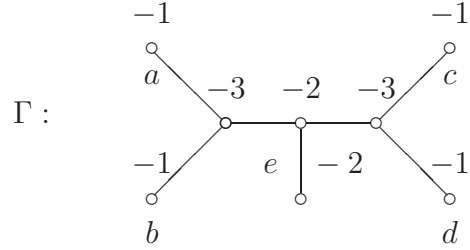
\vspace*{5mm}
$$
 \Gamma:\quad\qquad
 \cshiftup{-1}{a}\cshiftdown{-1}
 {b}\hspace{9.8mm}
 \nwlin\cou{-3}{}\swlin\co{}
\lin\co{-2}\slin \cshiftdown{\phantom{-}e\quad-2}{}
\lin\co{-3} \solin \nolin\hspace{9.8mm}
\cshiftup{-1}{c}\cshiftdown{-1}{d}
\quad
$$
\vspace*{10mm}

\noindent

\noindent admits two different contractions $\G\to A$ and $\G\to
B$ to the linear chains $A=B:=[[0,0]]$, namely by contracting the
subgraphs $\G\ominus \{a,b\}$, $\G\ominus\{c,d\}$, respectively.
We note that the extremal linear branch $T$ of $\Gamma$ consisting
just of the vertex $e$ is contracted in both $A$ and $B$, although
it is not contractible and $e\notin \hat A\cup\hat B$. Thus for a
pair of dominations $\Gamma\to A, \Gamma\to B$, an extremal linear
branch $T$ of $\Gamma$ contracted in both $A$ and $B$ is  not
necessarily contractible. \eexas

In the special case, where $\Gamma$ dominates two
circular graphs,
we have the following simple observation.

\blem\label{brpts1} If $\Gamma\to A$, $\Gamma\to B$ is a
relatively minimal pair of dominations of minimal circular graphs
$A,B$ then $\Gamma$ is as well circular. \elem

\bproof Since blowups do not change the topological type of the
graph, there is a unique circular subgraph
$\Gamma'\subseteq\Gamma$ dominating both $A$ and $B$. If
$\Gamma'\neq\Gamma$ then there is a branching point $c\in
B(\Gamma)$ on $\Gamma'$ and a branch $T$ at $c$ which is a
nonempty tree disjoint to $\Gamma'$. But then $T$ is contractible
and is contracted in both $A$ and $B$, which contradicts our
assumption of relative minimality. \eproof

\subsection{Admissible transformations}
Let us introduce the following notion.

\bdefi\label{creconstruction} A birational transformation of
weighted graphs $\Gamma\dasharrow\Gamma'$ which consists in a
sequence \bdi \gamma:\quad \Gamma=\Gamma_0&\rDashto^{\gamma_1}&
\Gamma_1 & \rDashto^{\gamma_2} &\cdots & \rDashto^{\gamma_n}&
\Gamma_n=\Gamma'\,\,, \edi where each $\gamma_i$ is either a
blowdown or a blowup, is called {\em admissible}\footnote{Or {\it
strict} in the terminology of \cite{Dai2}.} if the total branching
number remains constant at every step. For instance, a blowup
$\G'\to\G$ is admissible if it is inner or performed in an end
vertex of $\Gamma$, and a blowdown is admissible if its inverse is
so. Clearly, a composition of admissible transformations is
admissible.

More restrictively, we call $\gamma$ {\em inner} if the $\gamma_i$
are either admissible blowdowns or inner blowups. Thus the inverse
$\gamma^{-1}$ is admissible but not necessarily inner, see also
Definition \ref{graph.5} below. \edefi

The following proposition gives a precision of Theorem 3.2 in
\cite{Dai2}, which says that any birational transformation of
minimal graphs can be replaced by an admissible one.

\bprop\label{graph.9} If $\Gamma$ dominates two minimal graphs
$\Gamma_1$ and $\Gamma_2$ then there exists an admissible
transformation of $\Gamma_1$ into $\Gamma_2$ such that every step
is dominated by $\Gamma$. In other words, there is a birational
transformation of $\Gamma_1$ into $\Gamma_2$ such that each step
is dominated by $\Gamma$ and the total branching number stays
constant. \eprop

\bproof We may assume that $\Gamma$ minimally dominates
$\Gamma_1$ and $\Gamma_2$ that is, none of the $(-1)$-vertices of
degree $\le 2$ in $\Gamma$ is contracted in both $\Gamma_1$ and
$\Gamma_2$. Let $b$ be a branching point in $\Gamma$ not
contained in $B(\Gamma_1)\hat{}=B(\Gamma_2)\hat{}$. By Lemma
\ref{branching}(a), for $i=1,2$ there is a branch $C_i$ of
$\Gamma$ at $b$ that is contractible inside $\Gamma$ and is
contracted in $\Gamma_i$. Since $\Gamma$ dominates $\Gamma_1$ and
$\Gamma_2$ relatively minimally we have $C_1\ne C_2$, so these
are disjoint. Letting $\Gamma/C_i$ be the result of contracting
$C_i$ inside $\Gamma$ we obtain the diagram
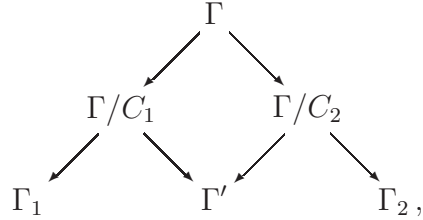
\bdi[size=1.5em]
   &&&&& \Gamma \\
   &&&&\ldTo &&\rdTo \\
   &&& \Gamma/C_1&&&&\Gamma/C_2 \\
   &&\ldTo&&\rdTo &&\ldTo &&\rdTo \\
   &\Gamma_1 &&&&\Gamma'&&& &\Gamma_2\,,
\edi where $\Gamma'$ is the minimal graph obtained from
$\Gamma/(C_1\cup C_2)$ by blowing down successively all
$(-1)$-vertices of degree $\le 2$. The total branching number of
$\Gamma/C_1$ and $\Gamma/C_2$ is strictly smaller than that of
$\Gamma$. Using induction on this number the result follows.
\eproof

\subsection{Elementary transformations} In this
subsection we assume that our weighted graph $\G$ is nonempty,
connected and has {\em only integral weights}.  In our principal
result (see Reduction Theorem \ref{graph.7} below) we use the
following operation on weighted graphs, cf. e.g. \cite{DaGi}.

\bdefi\label{graph.5} Given an at most linear vertex $v$ of $\G$
with weight $0$  one can perform the following transformations. If
$v$ is linear with neighbors $v_1$, $v_2$ then we blow up the
edge connecting $v$ and $v_1$ in $\G$ and blow down the proper
transform of $v$:  \begin{equation}\label{etr1} \bdi
\raisebox{1mm}{$\ldots
\cou{v_1}{w_1-1}\lin\cou{v'}{0}\lin\cou{v_2}{w_2+1}\ldots$ }\quad
&\rDashto & \quad \raisebox{1mm}{$\ldots
\cou{v_1}{w_1-1}\lin\cou{v'}{-1}\lin\cou{v}{-1}\lin\cou{v_2}{w_2}
\ldots $}\quad &\rTo& \quad \raisebox{1mm}{$
\ldots\cou{v_1}{w_1}\lin\cou{v}{0}\lin\cou{v_2}{w_2}\ldots
$}\quad. \edi \end{equation} Similarly, if $v$ is an end vertex
of $\G$ connected to the vertex $v_1$ then one proceeds as
follows: \be\label{etr2} \bdi \raisebox{1mm}{$\ldots
\cou{v_1}{w_1-1}\llin\cou{v'}{0} $} \quad &\rDashto&\quad
\raisebox{1mm}{$\ldots \cou{v_1}{w_1-1}\llin\cou{v'}{-1}\lin
\cou{v}{-1} $} \quad &\rTo &\quad \raisebox{1mm}{$\ldots
\cou{v_1}{w_1}\lin\cou{v}{0} $} \edi \quad.\ee These operations
(\ref{etr1}) and (\ref{etr2}) and their inverses will be called
{\em elementary transformations} of $\G$. If such an elementary
transformation involves only an inner blowup then we call it {\em
inner}. Thus (\ref{etr1}) and (\ref{etr2}) are inner whereas the
inverse of (\ref{etr2}) is not as it involves an outer blowup.
Clearly, elementary transformations are admissible in the sense
of Definition \ref{creconstruction}. \edefi

\bexas\label{graph.51} 1. If $[[w_0,\ldots, w_n]]$ is a weighted
linear graph with $w_k=0$ for some $k$, $0\le k\le n$ (see
\ref{birtr}) then \bdi[]
[[w_0,\ldots,w_{k-1}+1,\,0,\,w_{k+1}-1,\ldots w_n]] &\rDashto &
[[w_0,\ldots,w_{k-1},\,0,\,w_{k+1},\ldots w_n]] \,, \edi as well
as its inverse are elementary transformations. It is inner unless
$k=n$, and the inverse elementary transformation is inner unless
$k=0$.

2. Iterating inner elementary transformations as in (1),
for $a\in\Z$ we can transform a linear subgraph $L=[[w_1,0,w_2]]$
of a segment of $\G$ into $[[w_1-a,0,w_2+a]]$ leaving
$\G\ominus L$ unchanged (see
\cite[L. 4.15]{Dai2}). In particular, $[[w_1,0,0]]$ can be
transformed into $[[0,0,w_1]]$. Thus we can "move" pairs of
vertices of weight $0$ within segments of $\G$ by means of a
sequence of inner elementary transformations. \eexas

To simplify notation we will write $[[0_{m}]]$, $[[(a)_{m}]]$ for
the linear chain of length $m$ with all weights equal to $0$ or
$a$, respectively. Using the preceding examples we easily deduce
the following facts (cf. \cite[Lemmas 4.15-4.16]{Dai2}).

\blem\label{graph.6} Below $L$ stands for a linear subchain of  a
segment $\Sigma$ of $\G$.
By a sequence of
inner elementary transformations of $\G$ we can transform
$$
L=[[0_{2k},w_1,\ldots, w_n]] \quad\mbox{into}\quad
L^*=[[w_1,\ldots, w_n,0_{2k}]],\quad \forall k,n\ge1;\leqno(a)
$$
$$
L=[[w_1,0_{2k+1},w_2]] \quad\mbox{into}\quad [[w_1-a,0_{2k+1}, w_2+a]],\quad \forall a\in\Z ,\quad\forall k\ge 0;\leqno (b)
$$
$$
L=[[0_{2k+1},w_0,\ldots,w_n]] \quad\mbox{into}\quad
  [[0_{2k+1}, w_0-a,w_1,\ldots,w_n]],
  \quad \forall a, k,n\ge 0\,.\leqno (c)
$$
These elementary transformations leave $\G\ominus L$ unchanged
except possibly for the case (c), where the vertex connected to
the leftmost vertex of $L$ will change its weight. Allowing as
well outer elementary transformations (c) holds for all $a\in \Z$.
\elem

\bproof (a) follows by applying Example \ref{graph.51}(2)
repeatedly. To prove (b), by Example \ref{graph.51}(2) we can
transform $[[w_1,0_{2k+1},w_2]]$ into
$[[w_1,0_{2k-1},-a,0,w_2+a]]$ by a sequence of inner elementary
transformations. Hence (b) follows by induction. Finally, a chain
$L$  as in (c) can be transformed into $[[0,w_0,0_{2k},
w_1,\ldots,w_n]]$.  Applying elementary transformations as in
Definition \ref{graph.5}(\ref{etr2}) repeatedly we can transform
the latter chain into $[[0,w_0-a,0_{2k},w_1,\ldots,w_n]]$ and then
into $[[0_{2k+1},w_0-a,w_1,\ldots,w_n]]$. Note that in the case
$a\ge 0$ only inner transformations are required, see Example
\ref{graph.51}(1). \eproof

\subsection{Standard and semistandard graphs}

\bdefi\label{standard} A non-circular weighted graph $\G$ will be
called {\em standard} if each of its nonempty segments $L$ is one
of the linear chains \be\label{standard1}
[[0_{2k},w_1,\ldots,w_n]] \quad\mbox{or}\quad [[0_{2k+1}]] \,, \ee
where $k,n\ge 0$ and $w_i\le -2\,\forall i$. Similarly, a circular
graph will be called {\em standard} if it is one of the graphs
\begin{equation}\label{cstandard}
((0_{2k}, w_1,\ldots, w_n)), \quad\quad ((0_l, w))
\quad\mbox{or}\quad ((0_{2k}, -1,-1)) ,
\end{equation}
where $k,l\ge 0$, $n>0$, $w\le 0$ and $w_1,\ldots,w_n\le -2$.
We note that for $w=0$ the second graph in
(\ref{cstandard}) becomes $((0_{l+1}))$.

A graph
$\G$ will be called {\em semistandard} if each of its nonempty segments is
either a standard circular graph or one of the linear chains
\begin{equation}
\label{sstandard}
[[0_{l},w_1,\ldots,w_n]]
\quad\mbox{or}\quad [[0_{l},w_1,\ldots,w_n,0]]
\end{equation}
where  $l,n\ge 0$
and $w_i\le -2$ for $1\le i\le n$. Thus every standard graph is also semistandard.

For the standard linear chain $L=[[0_{2k},w_1,\ldots,w_n]]$ in
(\ref{standard1}) let us call the chain $L^*=[[0_{2k}, w_n,\ldots,
w_1]]=[[ w_1,\ldots, w_n,0_{2k}]]$ the {\em reversion} of $L$. By
Lemma \ref{graph.6}.a the reversion $L^*$ can be obtained from
$L$ by a sequence of inner elementary transformations.
\edefi

\brems\label{sstos} 1. {\em Every semistandard linear chain
can be transformed into a standard one by a
sequence of elementary transformations.} For it is one of the
chains in (\ref{sstandard}), and by Lemma \ref{graph.6} it can be
transformed into a chain as in (\ref{standard1}).

2. If $\Gamma=[[w]]$ ($w\neq -1$) has only one vertex then
it is either standard
or $w>0$. In the latter case after $w$ blowups we obtain the chain
$[[0,-1,-2,\ldots,-2]]$ of length $w+1$, which transforms as
before into a standard chain $[[0,0,-2,\ldots,-2]]$  of length $w+1$.
This transformation is not composed of elementary ones,
since it does not preserve the length.

3. Omitting an end vertex from a semistandard linear chain yields
again a semistandard chain. \erems

By a result of Daigle \cite[Thm. 4.23]{Dai2} every weighted tree
can be transformed into a unique standard one\footnote{A {\it
canonical} linear chain in the terminology of \cite{Dai2} is
different from our standard one, but in a recent version of
\cite{Dai2} these coincide.}, see also \cite{DaGi} for related
results. We give a simplified proof of this reduction result for
a general weighted graph in the following more precise form.

\bthm\label{graph.7} For a minimal graph $\Gamma$ the following
hold.
\begin{enumerate}
[(a)] \item If $\Gamma$ is not reduced to a point then $\G$ admits
an inner transformation into a semistandard graph, i.e.\ it can be
transformed into a semistandard one via a sequence of admissible
blowdowns and inner blowups. \item $\G$ allows an admissible
transformation into a standard graph.
\end{enumerate}\ethm

\bproof (b) follows from (a) in view of Remarks
\ref{sstos}.1 and 2. To show (a), after performing a suitable
sequence of inner blowups $\tilde\G\to\G$ we can achieve that all weights of vertices
inside the segments
of $\tilde\G$ become $\le 0$, and because of the minimality
assumption, these segments remain non-contractible. The result is
now a consequence of the following claim applied to each segment
$\Sigma$ of $\tilde \G$.

\smallskip\noindent
{\it Claim: A non-contractible segment $\Sigma $ with all weights
$\le 0$ can be transformed into one of the graphs in
(\ref{standard1}), (\ref{cstandard}) or (\ref{sstandard}) by
means of a sequence of admissible blowdowns and inner elementary
transformations.}\smallskip

\noindent {\it Proof of the claim.} We proceed by induction on the
length of $\Sigma$. Our assertion trivially holds if $\Sigma$ is
of length $1$. If $\Sigma$ contains a subgraph $[[w_1,-1,w_2]]$
with $w_1,w_2< 0$ then we can contract the $(-1)$-vertex to
obtain a segment with a smaller number of vertices.
Similarly, if $\Sigma=[[-1,w_1,\ldots]]$ with $w_1<0$ then again
we can contract the $(-1)$-vertex and get a segment of
smaller length.

Using Lemma \ref{graph.6} and contractions as above we can
collect inner vertices of weight 0 in $\Sigma $. If $\Sigma $ is
non-circular then we can also assume that they are on the left. In
this way we obtain one of the following nonempty graphs:
\begin{equation}\label{pp}
[[0_m,w_0,\ldots,w_n]], \quad\quad [[0_m,w_0,\ldots,w_{n+1},0]]
\quad\mbox{or}\quad ((0_m,w_0,\ldots,w_{n+1}))\,,
\end{equation}
where $m\ge 0$,
the sequence $(w_i)$ can be empty, $w_0, w_{n+1}\le -1$ and
$w_1,\ldots, w_n\le -2$.


Let us consider the first graph in (\ref{pp}). If $w_0\le -2$ then
it is semistandard. If $m=2k+1$ is odd  and $w_0=-1$ then by Lemma
\ref{graph.6}.c it can be transformed into
$[[0_{2k+1},-2,w_1,\ldots,w_n]]$ by a sequence of inner elementary
transformations. In the case $m=0$  and $w_0=-1$ we can contract
the $(-1)$-vertex and apply the induction hypothesis. If $m=2k>0$
is even and $w_1=-1$, then  we can first contract the
$(-1)$-vertex to obtain the chain
$[[0_{2k-1},1,w_1+1,\ldots,w_n]]$ and then transform this by inner
elementary transformations into $[[0_{2k},w_1+1,\ldots,w_n]]$, see
again Lemma \ref{graph.6}.c. Now the result follows by induction.

Similarly, combining the same operations and reversing the
ordering, if necessary, provides a reduction to a semistandard
graph in the other two cases in (\ref{pp}). The details are left
to the reader. \eproof

\subsection{Zigzags and standard zigzags} Let $V$ be a normal surface
and $X$ be a completion of $V$ by a divisor $D$ with simple normal crossings
(or by an SNC-divisor, for short) so that $D$ is contained in
the regular part $X_{\rm reg}$ of $X$.
The dual graphs $\Gamma (D)$ of such $D$
are restricted by the Hodge index theorem. We use the following
terminology.

\bdefi An SNC divisor $D\subseteq X_{\rm reg}$ with
irreducible components $C_1,\ldots, C_n$ in a complete
algebraic surface $X$ will be called a {\em zigzag} if the following conditions
are satisfied.

$\diamond$ The curves $C_i$ are rational $\forall i=1,\ldots,n$.

$\diamond$ The dual graph of $D$ is a linear chain
$L=[[w_0,w_1,\ldots, w_n]]$ such that the adjacency
matrix\footnote{See the Appendix below.} $I(L)$ has at most one
positive eigenvalue.

A zigzag will be called (semi)standard if its dual graph
(also called a zigzag)
has the corresponding property, see Definition \ref{standard}.
\edefi

We remind the reader that the number of positive eigenvalues is a
birational invariant of a graph, see \ref{adjma} in the Appendix
below. Thus a chain birationally equivalent to a zigzag is again a
zigzag. We also note that our terminology is different from the
one introduced in \cite{DaGi}. Indeed, in \cite{DaGi} by a
standard zigzag the authors mean an $m$-standard, quasistandard or
semistandard zigzag, whereas our standard zigzags are $0$-standard
in the sense of \cite{DaGi}. In the following lemma we describe
the dual graphs of (semi)standard zigzags (cf. \cite[Prop.
7.8]{Dai2}).

\blem\label{dgsz}
The possible standard zigzags are the chains
\be
\label{standardzigzag}&& [[0]]\,, \quad [[0,0,0]]
\,\quad\mbox{and}\quad
  [[0,0,w_1,\ldots,w_n]],\quad \mbox{where}
\quad n\ge 0,\, w_j\le-2\,\,\,\forall j, \ee whereas for the
semistandard ones we have additionally the possibilities
\be\label{mstandard} && [[0,w_1,\ldots,w_n]]\,, \,\,\,
[[0,w_1,0]], \qquad\mbox{where }
n\ge 0\quad\mbox{and}\quad w_j\le-2\,\,\,\forall j\,.\ee
\elem

\bproof This is an immediate consequence of the following
claim\footnote{Alternatively, one can derive the result from
Lemma \ref{cfr1}, cf. also Proposition \ref{numeig}.b.}.
\smallskip

\noindent {\em Claim: A zigzag $\G$ cannot contain two vertices
$v_i$ and $v_j$ of weight $\ge 0$ unless they are joined by an
edge, or $\G$ has at most length $3$. }
\smallskip

Otherwise, using the labelling as in (\ref{birtr}), we first
perform suitable inner blowups in $\Gamma$ to make the weights
$v_i^2=v_j^2=0$ and then, using Lemma \ref{graph.6},  we move the
$0$-weight on the left to $v_0$ and the one on the right to $v_n$.
By means of further elementary transformations we can assign to
$v_1$ and $v_{n-1}$ arbitrary weights e.g., $\ge 2$. But then the
adjacency matrix $I(\Gamma)$ has the symmetric submatrix
$\footnotesize
\begin{pmatrix} v_1^2 & v_1.v_{n-1}\\ v_1.v_{n-1} &
v_{n-1}^2\end{pmatrix}$, which is positive definite since
$v_1.v_{n-1}\le 1$ and $v_1^2, v_{n-1}^2\ge 2$. This contradicts
the assumption that $I(\Gamma)$ has at most one positive
eigenvalue, proving the claim. \eproof

\section{Birational transformations of standard graphs}

The central result of this section is the following structure
theorem for birational transformations of standard graphs.

\bthm\label{graph.8} Any birational transformation of one
semistandard graph into another one can be decomposed into a
sequence of elementary transformations.

More precisely, if $\Gamma\to A$, $\Gamma\to B$ is a pair of
dominations of semistandard graphs then $B$ can be obtained from
$A$ by a sequence of elementary transformations such that every
step is dominated by some inner blowup of $\Gamma$.\ethm

In Subsections 3.1 and  3.3 we completely describe admissible
transformations of linear and circular graphs, respectively. In
Subsection 3.4 the proof of \ref{graph.8} is then reduced to these
special cases.

\subsection{Admissible transformations of standard linear chains}
\bsit\label{dgr} In this subsection we consider a diagram
\be\label{diag1}
\bdi[size=1.5em]
&&& \Gamma &&&\\
&&\ldTo^{} &&\rdTo^{}&& \\
& A& & \rDashto & &B& \edi \ee where $\Gamma,A,B$ are linear
chains. We let "$<$" be an ordering on $\Gamma$, and we consider
the induced ordering of $A$ and $B$, respectively.  The main
results are Propositions \ref{p1} and \ref{p2} below, where we
describe completely all such birational transformations. One of
the key observations in the proofs is provided by the following
lemma. \esit

\blem\label{l1} Let $\Gamma$ be a linear chain, and let $\Gamma\to
A$, $\Gamma\to B$ be a pair of dominations, where the vertices
$a_1<a_2<\ldots <a_n$ and $b_1<b_2<\ldots <b_m$ of $A$ and $B$,
respectively, are ordered upon an ordering in $\G$. Then the
following hold.

(a) If $n\ge 2$ and $a_1^2=b_1^2=0$ then $m\ge 2$ and $\hat
a_2=\hat b_2$ in $\Gamma$.

(b) If $n,m\ge 3$ and for some $k,l$ with $1<k<n,\,\,\,1<l<m$ we
have $a_k^2=b^2_l=0$ then $a_{k-1}<b_{l-1}$ if and only if
$a_{k+1}<b_{l+1}$. \elem

\bproof To show (a), let us note first that in the case $m=1$ the
subchain $\Gamma_{<\hat a_2}$ is properly contained in $\G$, so by
Zariski's Lemma \ref{zar} its intersection matrix is negative
definite. As it contracts to $[[0]]$ in $A$ this gives a
contradiction. Hence $m\ge 2$.

To show the remaining assertion, let us suppose e.g.\ that $\hat
a_2<\hat b_2$. Then the subchain $\Gamma_{<\hat b_2}$ of $\Gamma$
is blown down to $[[0]]$ in $B$. By Zariski's Lemma \ref{zar},
every proper subgraph of $\Gamma_{<\hat b_2}$ is negative
definite. However, by assumption $\Gamma_{<\hat b_2}$ properly
contains $\Gamma_{<\hat a_2}$, and $\Gamma_{<\hat a_2}$ is not
negative definite as it is contracted to $[[0]]$ in $A$. This
contradiction proves that indeed $\hat a_2=\hat b_2$.

To deduce (b), we consider the open intervals $\G_A$ between
$\hat a_{k-1}$ and $\hat a_{k+1}$ and $\G_B$ between $\hat
b_{l-1}$ and $\hat b_{l+1}$ in $\G$. Clearly $\G_A$ and $\G_B$
are contracted to $[[0]]$ in $A$, $B$, respectively. If e.g. $\hat
a_{k-1}<\hat b_{l-1}<\hat b_{l+1}\le \hat a_{k+1}$ then the
inclusion $\G_B \subseteq \G_A$ would be proper contradicting
Zariski's Lemma \ref{zar}. The proof in the other cases is
similar. \eproof

The following results show that a non-trivial admissible
birational transformation between two standard linear chains can
exist exclusively for two chains of the same odd length with all
zero weights.

\bprop\label{p1} Let $\G,A,B$ be linear chains, and let $\G\to A$,
$\G\to B$ be a relatively minimal pair of dominations. Assume that
$A,B$ are standard and at least one of them is different from
$[[0_{2k+1}]]$ for all $k\ge 0$. If the groups of zeros in $A,B$
are both on the left then $\Gamma\to A$ and $\Gamma\to B$ are
isomorphisms. \eprop

\bproof
We denote as before by $a_1<a_2<\ldots <a_n$ and $b_1<b_2<\ldots
<b_m$  the vertices of $A$, $B$, respectively, ordered upon an
ordering in $\G$. We proceed by induction on the length of $A$.

If $a_i^2\le-2$ $\forall i$ then also $\hat a_i^2\le -2$ $\forall i$,
so  by relative minimality there is no
$(-1)$-vertex in $\Gamma$ blown down in $B$. Hence $\Gamma\to B$
is an isomorphism. Since $B$
is minimal, $\Gamma$ is minimal too and so $\Gamma\to A$ as well
is an isomorphism. The same conclusion holds in the case where $b_j^2\le-2$
$\forall j$.

Thus we may restrict to the case where $a_1^2=0=b_1^2$. By our
assumptions one of the chains $A$, $B$ has length $\ge 2$ and so,
by Lemma \ref{l1}(a), $n,m\ge 2$ and $\hat a_2=\hat b_2$ in
$\Gamma$. The graphs
$$
A'=A\ominus\{a_1,a_2\},\quad B'=B\ominus\{b_1,b_2\}\quad
\mbox{and}\quad
\Gamma'=\Gamma_{>\hat a_2}
$$
are linear, $A',B'$ are still standard and
$\Gamma'\to A'$, $\Gamma'\to B'$ is a relatively minimal
pair of dominations. Note that by our assumption, $A'\neq [[0_{2k-1}]]$ or
$B'\neq [[0_{2k-1}]]$. Using induction we get that $\Gamma'\to A'$,
$\Gamma'\to B'$
are isomorphisms. Now
$$
A''=A\ominus A',\quad B''=B\ominus B'
$$
are both dominated by $\Gamma_{\le \hat a_2}$. Since
$a_2^2=b_2^2=0$, by Lemma \ref{l1} $\hat a_1=\hat b_1$. Taking
into account the equality $\hat a_2=\hat b_2$ it follows that
$\hat a_1$ and $\hat a_2$ are neighbors in $\G$ and so $\Gamma\to
A$ and $\Gamma\to B$ are isomorphisms as required. \eproof

Next we describe the birational transformations of the standard
graphs $[[0_{2k+1}]]$.

\bdefi\label{move} For $A=[[0_{2k+1}]]$ with vertices
$a_1,a_2,\ldots,a_{2k+1}$, we let $\tau:A\dasharrow B$ be the
birational transformation consisting of the outer blowup at $a_1$
and the inner blowups at the edges $[a_{2i},a_{2i+1}]$,
$i=1,\ldots,k$ followed by the contraction of the vertices $\hat
a_{2i+1}$, $i=0,\ldots,k$. Thus $\tau$ fits the diagram
(\ref{diag1}) with $\Gamma=[[(-1)_{3k+2}]]$ and $B=[[0_{2k+1}]]$.
We call $\tau$ the {\it left move} and $\tau^{-1}$ the {\it right
move}. \edefi

A move $\tau$ admits a decomposition into a sequence of
elementary transformations consisting in the above blowups and
blowdowns, once at time.

\bexa\label{5vert} The linear chain
$$\G=[[(-1)_5]]\,:\qquad\cou{-1}{\hat b_1}\lin\cou{-1}{\hat a_1}
\lin\cou{-1}{\hat a_2=\hat b_2}\lin\cou{-1}{\hat
b_3}\lin\cou{-1}{\hat a_3}$$
\smallskip

\noindent dominates two linear chains
$$A=[[0_3]]\,:\qquad\cou{0}{a_1}\lin\cou{0}{a_2}
\lin\cou{0}{a_3}\qquad\mbox{and}\qquad
B=[[0_3]]\,:\qquad\cou{0}{b_1}\lin\cou{0}{b_2} \lin\cou{0}{b_3}$$

\noindent resulting in a left move $\tau: A \dashrightarrow B$.
\eexa

\bprop\label{p2} Let $A=[[0_{2k+1}]]$, $B=[[0_{2l+1}]]$ and let
$\Gamma$ be a linear chain. If $\Gamma\to A$, $\Gamma\to B$ is a
pair of dominations then $k=l$ and the resulting birational
transformation $A\dasharrow B$ is equal to $\tau^s$ for some
$s\in\Z$. In particular it is a composition of elementary
transformations. \eprop

\bproof We may assume that $k\le l$. With the notation as before,
by Lemma \ref{l1}(a) we have $\hat a_{2i}=\hat b_{2i}$ in $\Gamma$
$\forall i=1,\ldots,k$ and, similarly, $\hat a_{2k-2i}=\hat
b_{2l-2i}$ $\forall i=0,\ldots,k-1$. This is only possible if
$k=l$.  If $\hat a_1=\hat b_1$ in $\Gamma$ then applying Lemma
\ref{l1}(b) repeatedly we obtain that $\hat a_{2i+1}=\hat
b_{2i+1}$ for $i=1, \ldots, k$ and so $\hat A=\hat B$.

Now assume that $\hat a_1<\hat b_1$. Using again Lemma \ref{l1}(b)
gives that $\hat a_{2i+1}<\hat b_{2i+1}$ for $i=1, \ldots, k$.
Hence to obtain $\Gamma$ from $B$ all the edges
$[b_{2i},b_{2i+1}]$ are blown up, and furthermore an outer blowup
is performed at the vertex $b_1$. Consequently the linear chain
$B'=\tau (B)$ is dominated as well by $\Gamma$. Moreover the
distance of $\hat a_1$ to $\hat B'$ in $\Gamma$ is strictly
smaller than the distance to $\hat B$. Using induction on this
number the assertion follows. \eproof

\subsection{Linear dominations of semistandard chains}
It follows from Propositions \ref{p1}, \ref{p2} and Remark
\ref{sstos} that every birational transformation $A\dashrightarrow
B$ of semistandard linear chains dominated by a linear chain $\G$
can be decomposed into a sequence of elementary transformations.
Later on we need the stronger result which says that one can
dominate $A, B$ and the intermediate graphs even by a suitable
{\em inner} blowup of $\Gamma$. To deduce this fact we need the
following observation.

\blem\label{anal0}
Let $\G'\to\G$ be a domination of weighted graphs and let
$\gamma: \tilde\G\dasharrow \G$ be an admissible transformation.
Then there is a commutative diagram
$$ \bdi[small]
\tilde \G'&\rTo^{\gamma'}& \G'&\\
\dTo && \dTo\\
\tilde\G&\rDashto^{\gamma}&\G\,, \edi
$$
where the solid arrows are dominations and $\gamma'$
is admissible. Moreover, if $\gamma$ is inner then also $\gamma'$
can be chosen to be inner.
\elem

\bproof Decomposing $\gamma$ into a sequence of admissible blowups
and blowdowns, it is enough to consider the following 3 cases:
(i)  $\gamma$ is an admissible blowup in $\tilde \G$ of a vertex
$v$ of $\G$, (ii) $\gamma$ is an inner blowup at an edge
$e=[v_1,v_2]$ of $\G$, or (iii) $\gamma$ is an outer blowup at an
end vertex $v_0$ of $\G$. So in case (i) $\tilde\G$ has 1 vertex
more than $\G$, and in cases (ii), (iii) it has 1 vertex less than
$\G$.

In case (i) $\G$ and then also $\G'$ dominates $\tilde \G$, so we
can choose $\tilde \G'=\G'$ and $\gamma'=$id. In case (ii), if the
edge $e$ is blown up in $\G'$ then we can again choose
$\tilde\G'=\G'$ and $\gamma'=$id. Otherwise the proper transforms
$\hat v_1$ and $\hat v_2$ of $v_1$, $v_2$ in $\G'$ are neighbors.
Blowing up the edge $[\hat v_1,\hat v_2]$ we obtain a graph
$\tilde\G'$ with the desired properties.

Similarly, if in case (iii) $\G'\to\G$ factors through an outer
blowup at $v_0$, we can choose $\tilde\G'=\G'$ and $\gamma'=$id.
Otherwise the proper transform $\hat v_0$ of $v_0$ in $\G'$
remains an end vertex, and we can choose $\gamma'$ to be the outer
blowup at $\hat v_0$. \eproof

\brem\label{rmi} The pair of dominations $\tilde \G'\to
\Gamma',\,\tilde \G'\to \tilde \G$ is not necessarily relatively
minimal. Indeed, the same vertex can appear or disappear under
$\G'\to \G$ and $\tilde \G\dashrightarrow \G$, cf. Remark
\ref{philo} below. \erem

Now we can deduce the following result.

\bprop\label{c1} If a linear chain $\Gamma$  dominates two
semistandard chains\footnote{See Definition \ref{standard}.} $A$
and $B$ then $B$ can be obtained from $A$ by a sequence of
elementary transformations such that every step is dominated by
some inner blowup of $\Gamma$. \eprop

\bproof As noted at the beginning of this subsection, $A$ can be
obtained from $B$ by a sequence of elementary transformations and
so, in particular, $A$ and $B$ have the same length. We denote the
vertices of $A$ and $B$ as before by $a_1<\ldots<a_n$ and
$b_1<\ldots<b_n$, respectively.

We proceed by induction on $n$. If both chains $A$ and $B$ are standard then the result is a consequence of Propositions \ref{p1} and \ref{p2}.
In particular this settles the case $n=1$.
Let us assume for the rest of the proof that $n\ge 2$.

If all weights in $A$ or in $B$ are $\le -2$ then
$\hat A=\hat B$ similarly as in the proof of Proposition
\ref{p1}. Hence we may assume in the sequel that both $A$ and $B$
have zero weights. Up to interchanging $A$ and $B$ or reversing both chains we have to consider the following two cases.

Case 1: $a_1^2=b_1^2=0$, \qquad Case 2: $a_1^2=b_n^2=0$, but $a_n^2\neq 0$ and $b_1^2\neq 0$.

In case (1), by Lemma \ref{l1}.a we have $\hat a_2=\hat b_2$ and
so, as in the proof of Proposition \ref{p1},
$\Gamma'=\Gamma_{>\hat a_2}$ dominates the semistandard chains
$$A'=A\ominus\{a_1,a_2\}\quad\mbox{and}\quad B'
=B\ominus\{b_1,b_2\}\,,$$ while $\Gamma''=\Gamma_{<\hat a_2}$
dominates the chains $A''=\{a_1\}$ and $B''=\{b_1\}$. Applying the
induction hypothesis and the case $n=1$, $B'$, $B''$ can be
obtained from $A'$, $A''$, respectively, by a sequence of
elementary transformations such that every step is dominated by
some inner blowup of $\Gamma'$, $\G''$, respectively. After
performing the same sequence of elementary transformations in $A$
we may assume that $A'=B'$ and $A''=B''$ and so we are done.

In case (2) we have
$$
A=[[0,\alpha,\ldots]]\quad\mbox{and}\quad B=[[\ldots,\beta,0]]\,.
$$

$\diamond$ If $\alpha=a_2^2=0$ then we can move this pair of zeros
in $A$ to the right by a sequence of inner elementary
transformations and so reduce to the case already treated. Indeed,
by Lemma \ref{anal0} the resulting new chain $A$ is again
dominated by an inner blowup of $\Gamma$. Thus we may suppose that
$\alpha=a_2^2\le -2$ and, symmetrically, $\beta=b_{n-1}^2\le -2$,
so $A$ and $B$ have the only zero weights $a_1^2=b_n^2=0$.

$\diamond$ If $\hat a_1$ is not an end vertex then $\G$ dominates
as well the chain $A_1=[[0,\alpha+1,\ldots]]$ obtained from $A$ by
an elementary transformation involving an outer blowup. Using
decreasing induction on $-\alpha$ the result can thus be reduced
either to the case treated above where $a_1^2=a_2^2=0$, or to the
case where $\hat a_1$ is the leftmost vertex of $\Gamma$. Clearly
the same argument applies to $b_n$ provided that $\hat b_n$ is not
an end vertex of $\G$.

 $\diamond$ So we are finally reduced to the case where
both $\hat a_1$ and $\hat b_n$ are end vertices of $\G$, $\hat
a_1$ is the leftmost one and $\hat b_n$ is the rightmost one. We
may suppose that the pair of dominations $\Gamma\to A,\,\Gamma\to
B$ is relatively minimal so that every $(-1)$-vertex in $\Gamma$
is contained in $\hat A\cup \hat B$. By virtue of
(\ref{sstandard}) $\hat a_1$ and $\hat b_n$ are the only possible
$(-1)$-vertices in $\Gamma$. If $\hat a_1^2=0$ then $\hat a_1$
cannot be contracted in $B$ and its image in $B$ is the only
$0$-vertex $b_n$. Therefore $\hat a_1=\hat b_n$ is the only end
vertex of $\Gamma$, so $A=B=[[0]]$, which contradicts our
assumption that $n\ge 2$. Thus $\hat a_1^2<0$, and passing from
$A$ to $\Gamma$ an inner blowup at the edge $[a_1,a_2]$ was
performed. Hence $\hat a_2^2<a_2^2=\alpha$ and, likewise, $\hat
b_n^2<b_n^2=\beta$, $a_1^*\neq\{\hat a_1\}$ and $b_n^*\neq\{\hat
b_n\}$, where $a^*$ stands for the total transform of $a$ in
$\Gamma$. Moreover $a_1^2=0$ forces that $\hat a_1$ must be
contained in $b_n^*$ and so $\hat b_{n-1}<\hat a_1$ in $\G$.
Similarly, $\hat b_n$ must be contained in $a_1^*$ and  so $\hat
b_n<\hat a_2$  in $\G$ contradicting the fact that $\hat a_1$ and
$\hat b_n$ are the end vertices of $\G$. This concludes the proof.
\eproof

\subsection{The circular case} In this subsection we treat
birational
transformations of standard circular graphs. Unlike in the linear
case, there are many birational transformations of such graphs as
soon as they contain $0$-vertices. Nevertheless, all these
transformations are composed of simple ones, which we call {\em
turns} and {\em shifts}, see Theorem \ref{maincirc} below. Let us
introduce the following notions.

\bdefi\label{tech} A circular graph $A$ will be called {\em almost
standard} if it is standard with all weights $<0$ \footnote{that
is (see (\ref{cstandard})) $A$ is one of the graphs $((w))$, $
w\le -1,$ $((-1,-1))$, $((\alpha_1,\ldots,\,\alpha_n))$,
$\alpha_i\le -2\,\forall i.$} or if it can be written in the form
$$
A=((0_{2k+1},\,\alpha_0,\ldots,\,\alpha_n)), \quad k\ge 0,
$$
where the sequence $(\alpha_0,\ldots,\,\alpha_n)$ can be either
empty, or equal to one of
$$(w)\quad\mbox{with}\quad w\le 0,\quad (0,-1),\quad (-1,0),\quad
(0,-1,-1)\quad\mbox{or}\quad (-1,-1,0)\,,$$ or satisfies the
conditions $$n\ge 1,\quad\alpha_0,\,\alpha_n\le 0,\quad
\alpha_0+\alpha_n\le -2\quad\mbox{and}\quad\alpha_i\le
-2\,\,\forall i=1,\ldots, n-1\,.$$ \edefi

\brems\label{remzeros} 1. We note that for any sequence
$(\alpha_0,\ldots,\alpha_n)$ as in Definition \ref{tech} the
standard form of the chain $[[\alpha_0,\ldots,\alpha_n]]$ contains
at most two zeros. Hence its intersection form has at most 1
positive eigenvalue, see Proposition \ref{numeig} in the Appendix.

2. Slightly more generally, the standard form of the chain
$[[0_{2l},\alpha_0,\ldots,\alpha_n]]$ contains at most $2l+2$
zeros. \erems

Clearly every standard circular graph is almost standard.  To
treat birational transformations of standard graphs it is
convenient to consider these more generally for almost standard
graphs. In a special case their classification is simple.

\blem\label{circl3} Let $A$, $B$ be almost standard circular
graphs and let $\Gamma\to A$ and $\Gamma\to B$ be dominations. If
all weights of $A$ or $B$ are $\le -1$, then $\hat A=\hat B$.
\elem

\bproof If all weights of $A$ or $B$ are $\le -2$ then with the
same argument as in the proof of Proposition \ref{p1} it follows
that $\hat A=\hat B$. This argument also works if $A$ or $B$ is
one of the circular graph $((-1,-1))$, $((-1))$, whence the
result. \eproof

In the case that $A$ or $B$ contain 0-vertices the situation is
much more complicated.  In analogy with moves of the linear chains
$[[0_{2k+1}]]$ (see Definition \ref{move}), we introduce the
following operation on such almost standard circular graphs.

\bdefi\label{defshift} Let
$A=((0_{2k+1},\alpha_0,\ldots,\alpha_n))$ be an almost standard
circular graph.
By a {\em shift} we mean  the birational transformation
$$
\sigma: A=((0_{2k+1},\alpha_0,\ldots,\alpha_n))\dasharrow A'
=((0_{2k+1},\,\alpha_0-1,\,\alpha_1,\,\ldots,\,\alpha_{n-1},\,
\alpha_n+1))
$$
composed of a sequence of elementary transformations, which send
the almost standard graph
$$
A=((0_{2k+1},\,\alpha_0,\ldots,\,\alpha_n))=
((\alpha_n,0_{2k+1},\,\alpha_0,\ldots,\,\alpha_{n-1}))
$$ as in \ref{tech} with $n\ge 1$
into
$$
((\alpha_n+1,0_{2k+1},\,\alpha_0-1,\,\alpha_1,\,\ldots,\,\alpha_{n-1}))
=((0_{2k+1},\,\alpha_0-1,\,\alpha_1,\,\ldots,\,\alpha_{n+1},\,\alpha_{n}+1))
\,,$$ see Lemma \ref{graph.6}(b).

The resulting graph $A'=\sigma(A)$ is almost standard provided
that  $\alpha_n\le -1$. Thus the shift of at least one of the
graphs $A$ or
$$
A^*:=((0_{2k+1},\,\alpha_n,\,\alpha_{n-1},\,\ldots,\,\alpha_1,\,\alpha_0))
$$
is again almost standard. Since $A$ can be written in the form
$A^*$, the inverse of a shift is again a shift. \edefi

\brems\label{remshift}  1. Any almost standard circular graph
$((0_{2k+1},\,\alpha_0,\ldots,\alpha_{n-1},\,\alpha_n))$ with
$n\ge 1$ can be transformed into a standard one
$((0_{2k+2},\,\alpha_1,\ldots,\alpha_{n-1},\,\alpha_n+\alpha_0))$
by a sequence of shifts.

2. A shift transforms a standard circular graph
$A=((0_{2k},\alpha_1,\ldots, \alpha_n))$ with $k\ge 1$ and
$\alpha_1,\ldots, \alpha_n\le -2$ into
$((0_{2k-1},-1,\alpha_1,\ldots, \alpha_{n-1},\alpha_n+1))$. Thus
by a sequence of shifts we can transform $A$ into the standard
graph
$$
((0_{2k-1},\alpha_n,\alpha_1\ldots,
\alpha_{n-1},0))=((0_{2k},\alpha_n,\alpha_1,\ldots,
\alpha_{n-1})).
$$
Hence, if for $A$ as above and $B=((0_{2k},
\beta_1\ldots,\beta_n))$, the sequences
$(\alpha_1,\ldots,\alpha_n)$ and $(\beta_1,\ldots,\beta_n)$ are
equal up to a cyclic permutation and reversion, then $A$ and $B$
are birationally equivalent via a sequence of shifts.

For instance, the standard circular graphs
$$((0,0,-3,-5,-2)),\quad ((0,0,-2,-3,-5))\quad\mbox{and}\quad
((0,0,-5,-2,-3))$$ obtained one from another by  cyclic
permutations of the nonzero weights, are birationally equivalent
via a sequence of shifts.

3. Let $a_1, \ldots, a_{2k+1}, a_{2k+2},\ldots a_{2k+n+2}$ be the
vertices of $A$ numbered according to the ordering of weights in
$((0_{2k+1},\alpha_0,\ldots,\alpha_n))$. Then the vertices $a_i$
with $i$ even or $i\ge 2k+2$ are not blown down by a shift
$\sigma$ as in Definition \ref{defshift}. More precisely, if $A$
and $A'$ are dominated by a graph $\G$, where $A'=\sigma
(A)=((0_{2k+1},\alpha_0-1,\ldots,\alpha_n+1))$ has vertices $a_1',
\ldots, a_{2k+1}', a'_{2k+2},\ldots a'_{2k+n+2}$, then $\hat
a_{i}=\hat a_i'$ in $\G$ for $i$ even or $i\ge 2k+2$. Indeed,
performing elementary transformations of a chain $[[w,0,w']]$ at
the $0$-vertex in the middle, the two outer vertices will not be
blown down.

4. Implicitly, the definition of shifts addresses as well cyclic
renumbering. However, ignoring this procedure will not create
serious ambiguity. \erems

\bdefi\label{defturn} Let
$A=((0_{2k+1},\alpha_0,\ldots,\alpha_n))$ be an almost standard
circular graph.
A {\em turn} $\tau: A\dasharrow AÇ$
 consists in a sequence
of elementary transformations sending $A$ first into
$$((0_{2k-1},\alpha_0,0,0,\alpha_1,\ldots,\alpha_n)),\qquad\mbox{then into}
\qquad
((0_{2k-1},\alpha_0,\alpha_1,0,0,\alpha_2,\ldots,\alpha_n))\,$$
(see Lemma \ref{graph.6}.a), until we arrive at
$$
A'=((0_{2k-1},\alpha_0,\ldots,\alpha_n,0,0))=
((0_{2k+1},\alpha_0,\ldots,\alpha_n))\cong A\,.
$$
The inverse birational transformation will also be called a turn.
\edefi

In analogy with Remark \ref{remshift}(2) we make the following
observation.

\brem\label{remturn} Let $a_1, \ldots, a_{2k+1}, a_{2k+2},\ldots
a_{2k+n+2}$ be the vertices of $A$ ordered correspondingly to the
weights $((0_{2k+1},\alpha_0,\ldots,\alpha_n))$. Assume that $\G$
dominates both $A$ and
$A'=((0_{2k-1},\alpha_0,\ldots,\alpha_n,0,0))$, where $A'$ with
vertices $a_1', \ldots, a'_{2k+n+2}$ is obtained from $A$ by a
turn as in Definition \ref{defturn}. Then in $\G$ we have $\hat
a_{i}=\hat a_i'$ for $i\le 2k$ and $i=2k+n+2$. \erem

The following theorem gives a complete description of birational
transformations between standard circular graphs.

\bthm\label{maincirc} Any birational transformation of standard
circular graphs $A\dasharrow B$ is either an isomorphism or it can
be written as a composition of turns and shifts. \ethm

Before turning to the proof we mention the following corollary.

\bcor\label{cornew}
Any birational transformation of standard
circular graphs $A\dasharrow B$ can be written as sequence of
elementary transformations.
More precisely, if $\G\to A$ and $\G\to B$ are dominations
then $A$ can be obtained from $B$ by a sequence of elementary
transformations such that every step is dominated by some inner
blowup of $\G$.
\ecor

\bproof The first part is an immediate consequence of Theorem
\ref{maincirc}. The second part follows from this in view of
Lemmata \ref{brpts1} and \ref{anal0}. \eproof

Theorem \ref{maincirc} is shown in \ref{circp} and is a
consequence of Lemmata \ref{circl3} above and \ref{circl2},
\ref{circl5} below. The strategy of the proof is as follows.
Consider a pair of dominations $\G\to A$ and $\G\to B$, where $A$
and $B$ are standard circular graphs. By Lemma \ref{circl2} below,
applying shifts we may achieve that $\hat A\cap\hat B\ne
\emptyset$, and moreover, $\hat a=\hat b$ for some vertices $a$ of
$A$ and $b$ of $B$ with $a^2=b^2=0$. In a second step we reduce
the statement to the linear case by restricting to $A\ominus
\{a\}$,  $B\ominus \{b\}$ and $\G\ominus \{\hat a=\hat b\}$.

The following simple example shows that the case $\hat A\cap \hat
B=\emptyset$ indeed occurs.

\bexa Contracting in $\Gamma=((-3,-1,-2,-2,-1))$ alternatively the
subchains $[[-3,-1,-2]]$ or $[[-2,-2,-1]]$ yields the standard
circular graphs $A=B=((0,0))$, and we have $\hat A\cap \hat
B=\emptyset$. \eexa

To deal with the case $\hat A\cap\hat B= \emptyset$ we introduce a
portion of notation.

\bsit\label{tech2} We let $\Gamma\to A,\,\Gamma\to B$ be a
relatively minimal pair of dominations of almost standard circular
graphs $A$ and $B$ satisfying  $\hat A\cap\hat B=\emptyset$. We
also let $\hat A_1,\ldots,\hat A_s$ and $\hat B_1,\ldots,\hat B_t$
be the connected components of the graphs $\hat A$ and $\hat B$,
respectively. Then $A_i$, $B_i$ are connected subchains of $A$,
$B$, respectively, so that
$$
A=((A_1,\ldots,A_s))\quad\mbox{and, similarly,}
\quad B=((B_1,\ldots,B_t))\,.
$$
The nonempty linear subchain, call it $X_i$, between $\hat A_i$
and $\hat A_{i+1}$ in $\Gamma$ is contracted in $A$, so it
contains a $(-1)$-vertex. By the relative minimality assumption
this vertex must be in $\hat B$. Hence $X_i$ includes some
component $\hat B_j$. Similarly, between $\hat B_j$ and $\hat
B_{j+1}$ there is a unique component $\hat A_i$. This implies that
$s=t$, and with an appropriate enumeration of the components $B_i$
we can write
$$
\Gamma=((E_1\,\hat A_1\,F_1\,\hat  B_1\,E_2\,
\hat A_2\,F_2\,\hat B_2\,\ldots\,
\, E_s\,\hat A_s\,F_s\,\hat
B_s))\,.
$$
Using indices in $\Z/s\Z$ (so $E_{i+s}=E_i\,\,\forall i$) the
chains
$$
X_i=F_i\hat B_iE_{i+1}\qquad\mbox{and}\qquad Y_i=E_i\hat A_iF_i
$$
are contractible, so they contain at least one $(-1)$-vertex lying
in $\hat B_i$ and $\hat A_i$, respectively. \esit

\blem\label{circl1} With the assumptions and notation as in
\ref{tech2} the following hold. \bnum \item $\hat A_i$ and $A_i$
have no vertex of weight $\ge 0$ in their interior. \item $\hat
A_i$ contains a $(-1)$-vertex but not the string $[[-1,-1]]$.
Consequently, $A_i$ has no chain $[[-1,-1]]$ in its interior.
\item Each $A_i$ has at least 2 vertices. \item It is not possible
that $\hat A_i=[[\ldots, -1]]$ and $\hat A_{i+1}=[[-1,\ldots]]$
simultaneously. \item One of the end vertices of $A_i$ has weight
0. \enum The same assertions hold for the components $B_i$. \elem

\bproof (1) and (2) follow from the fact that $\hat A_i$ is part
of a contractible chain and so every subchain of $\hat A_i$ has a
negative definite intersection matrix.

To show (3), assume that $A_i$  has only one vertex. As $\hat
A_i$ contains a $(- 1)$-vertex, necessarily $A_i=[[0]]$ in this
case. But then to obtain $\G$ from $A$ a blowup would occur only
on one side of $A_i$, which leads to a contradiction since the
chains $X_i=F_j\hat B_j E_{j+1}$ are all nonempty.

If (4) were violated then between $A_i$ and $A_{i+1}$ just one
blowup would occur, so $B_i$ would be a chain of length 1
contradicting (3).

To deduce (5), if for some $i$ none of the end vertices of $A_i$
were of weight 0 then by (1) the end vertices of $\hat A_i$ would
be of weight $\le -2$. Hence a $(-1)$-vertex of $\hat A_i$ must
lie in the interior of $\hat A_i$ and so it is a $(-1)$-vertex in
$A_i$ too. Thus $A=((0_{2k+2},-1,-1))$ or
$A=((0_{2k+1},\alpha_0,\alpha_1,\ldots,\alpha_n))$ with
$\alpha_0=-1$ or $\alpha_n=-1$ in Definition \ref{tech}. However
in both cases any connected linear subchain of $A$ satisfying (1)
with a $(-1)$-vertex in its interior has an end $0$-vertex.
\eproof

\blem\label{circl2} Let $A$ and $B$ be almost standard circular
graphs, and let $\Gamma \to A$, $\Gamma\to B$ be a relatively
minimal pair of dominations with $\hat A\cap \hat B=\emptyset$.
Then $A$ can be transformed by a finite sequence of shifts all
dominated by $\G$ into a new almost standard circular graph $A'$
satisfying $\hat A'\cap \hat B\ne\emptyset$. \elem

\bproof By Lemma \ref{circl3} $A$ is almost standard of the form
$A =((0_{2k+1},\alpha_0,\ldots,\alpha_n))$ as in Definition
\ref{tech}. For $n\ge 1$, since $\alpha_0+\alpha_n\le -1$ and $A$
can be written in the form
$((0_{2k+1},\alpha_n,\ldots,\alpha_0))$, we may assume that
$\alpha_n\le -1$. Using Lemma \ref{circl1}(1) and (3), at least
$k$ among the components $A_i$, say, $A_1,\ldots,A_k$, are equal
to $[[0,0]]$. If $\alpha_0\le -1$ then by virtue of Lemma
\ref{circl1}(5), up to reversion there is just one extra component
$A_{k+1}=[[\alpha_0,\ldots,\alpha_n,0]]$. If $\alpha_0=0$ then
either there is just one extra component
$A_{k+1}=[[0,\alpha_1,\ldots,\alpha_n,0]]$, or $n\ge 2$ and there
are 2 extra components
$$
A_{k+1}=[[0,\alpha_1,\ldots,\alpha_l]]\quad\mbox{and}\quad
A_{k+2}=[[\alpha_{l+1},\ldots,\alpha_n,0]],\quad\mbox{where}\quad1\le
l\le n-1\,.
$$
Actually the latter case cannot occur. Indeed by Lemma
\ref{circl1}(4), if one of the chains $\hat {A_i}$ ends with a
$(-1)$-vertex, say on the right, then  all of them, in particular
$A_{k+1}$, will have a $(-1)$-vertex on the right. But this is
impossible since $\alpha_i\le -1\,\,\forall i=1,\ldots,n$.

Thus for any $n\ge 0$,
$$
A_i=[[0,0]]\quad\forall
i=1,\ldots,k\qquad\mbox{and}\qquad
A_{k+1}=[[\alpha_0,\,\ldots,\,\alpha_n,0]]\,.
$$
Let us first assume that $\hat A_i$ has a vertex of weight $-1$ on
the right for at least one $i$ in the range $1\le i\le k+1$.
Clearly $\G$ dominates the graph obtained from $A$ by blowing up
all edges between $A_i$ and $A_{i+1}$ (using cyclic indices as
before). Therefore $\G$ dominates both $A$ and the semistandard
circular graph
$$
A':=\sigma (A)=((0_{2k+1},\,\alpha_0-1,
\,\alpha_1,\,\ldots,\,\alpha_{n-1},\,\alpha_n+1))
$$
obtained from $A$ by a shift.
Consider the blowdown $\G'$ of the $(-1)$-vertex on
the right of $\hat A_i$. By
construction $A'$ and $B$ are still dominated by $\G'$. Applying
induction on the number of vertices of $\G$, the result follows.

For the rest of the proof we may, and we will, assume that $\hat
A_i$ has no $(-1)$-vertex on the right, so $\hat A_i=[[-1, w_i]]$
for $i=1,\ldots, k$. If $\alpha_0=0$, then by symmetry the same
argument as before works. Thus we may assume that $\alpha_0<0$.
Using Lemma \ref{circl1}(4) $\hat A_{k+1}$ has no end vertex of
weight $-1$. Since it contains a $(-1)$-vertex (see Lemma
\ref{circl1}(2)) it follows that this vertex lies in the interior
of $\hat A_{k+1}$ and so it is also a $(-1)$-vertex of $A_{k+1}$.
As $\alpha_n\le -1$ and inspecting \ref{tech} this is only
possible when $\alpha_n=-1$. Now we can reduce to the case already
treated: $\G$ dominates both $A$ and the semistandard graph
$A'=((0_{2k+1},\,\alpha_0-1,
\,\alpha_1,\,\ldots,\,\alpha_{n-1},0))=((0_{2k+1}, \alpha_0',
\ldots, \alpha_n'))$, where
$$
(\alpha_0', \ldots, \alpha_n')=(0,\alpha_0-1,\alpha_1,
\ldots, \alpha_{n-1})\,.
$$
Replacing $A$ by $A'$ and applying the previous case the result
follows. \eproof

\bsit\label{tech3}
Because of Lemma \ref{circl3} we may, and we will,
assume in the rest of the proof of Theorem \ref{maincirc}
that $A$ and $B$ are both almost standard with a vertex of weight 0. Let $\Gamma\to A$ and $\Gamma\to B$ be a pair of dominations with $\hat A \cap \hat B\ne \emptyset$. We fix an orientation of $\G$ and orient $A$
and $B$ accordingly. With these orientations let us write
$$
A=((0_{2k+1},\alpha_0,\ldots, \alpha_n))
\quad\mbox{and}\quad
B=((0_{2l+1},\beta_0,\ldots, \beta_m))
$$
as in \ref{tech}. Using the same ordering let us denote the
vertices of
$A$ and $B$ by $a_1, \ldots, a_{2k+n+2}$ and $b_1,\ldots,
b_{2l+m+2}$, respectively.
\esit

With this notation we have the following result.

\blem\label{lemshift}
For any vertex $a$ of an almost standard circular graph
$$
A=((0_{2k+1}, \alpha_0,\ldots, \alpha_n))\quad\mbox{with}\quad k,n\ge 0,
$$
there exists a sequence of shifts and turns transforming $A$ into
an almost standard circular graph $A'=((0_{2k+1},
\alpha_0',\ldots, \alpha_n'))$ so that at any step $a$ is not
blown down and one of the following conditions is satisfied.
\bnum[(i)] \item $\alpha_n'\le -1$ and the image of $a$ in $A'$ is
the 0-vertex $a_1'$; or

\item $n=0$ and the image of $a$ in $A'$ is the vertex $a'_{2k+2}$ of weight $\alpha_0'$.
\enum

\noindent
Moreover, if $A$ was standard then $A'$ can be assumed to be standard.
\elem

\bproof If $n\ge 1$ then we may suppose that $\alpha_n\le -1$,
since otherwise we can write $A$ in the form
$((0_{2k+1},0,\alpha_0,\ldots, \alpha_{n-1}))$. Let $a=a_i$ with
the notations as in  \ref{tech3}. According to the value of $i$ we
consider the following cases (a)-(d).

(a)  $i= 2k+2+j$, where $0\le j\le n$. If $n=0$ then $j=0$ and
$a=a_{2k+2}$ as needed in (ii). Otherwise according to Remark
\ref{remshift}(1) and (3) there is a sequence of shifts which do
not contract $a_i$ and transform $A$ into a standard circular
graph
$$
A_1=((0_{2k+1},\,
\alpha_0+\alpha_n,\,\alpha_1,\ldots,\alpha_{n-1},0))
=((0_{2k+1},0,\alpha_0+\alpha_n,\,\alpha_1,\ldots,\alpha_{n-1}))\,.
$$
If $i=2k+n+2$ then the image of $a$ in $A_1$ occupies the 1-st
position, as required in (i). Otherwise we continue as in Remark
\ref{remshift}(2) and transform $A_1$ by a sequence of shifts into
a standard graph
$$
A'=((0_{2k+2},\alpha_{j},\ldots,\alpha_{n-1},\,
\alpha_0+\alpha_n,\,\alpha_1, \ldots,\alpha_{j-1}))
$$
moving the vertex $a$ to the 1-st position.

(b) Suppose now that  $1\le i=2l+1\le 2k+1$.
By a sequence of turns
moving the $2l$ zeros on the left of $a_i$ to the right,
we can transform $A$
into
$$A_1=((0=a^2,\,0_{2k},\,\alpha_0,\,\alpha_1,\ldots,\alpha_{n}))$$
so that $a$ is not blown down and is placed onto the 1-st
position, as needed.

(c) Suppose that $n\ge 1$ and $1\le i=2l\le 2k+1$.  In view of
Remark \ref{remshift}(3), $a$ is not blown down under the shifts
which transform $A$ into a standard circular graph
$$A_1=((0_{2k+1},\,
\alpha_0+\alpha_n,\,\alpha_1,\ldots,\alpha_{n-1},0))=
((0_{2k+2},\,\alpha_0+\alpha_n,\,\alpha_1,\ldots,
\alpha_{n-1}))\,.$$ Now the vertex $a=a_i=a_{2l}$ has been moved
to the position $2l+1$. This provides a reduction to the previous
case, so we can further move $a$ to the 1-st position.

(d) Finally, if $n=0$ and $1\le i=2l\le 2k+1$ then $a$ can be
moved to the last position by a sequence of $k-l+1$ turns, see
Remark \ref{remturn}.

The reasoning in all 4 cases shows that if we start with a
standard circular graph $A$ then also $A'$ will be standard,
finishing the proof of the lemma.
\eproof

We are now ready to complete the proof of Theorem \ref{maincirc}
by the following lemma.

\blem\label{circl5}
With the notations and assumptions as in \ref{tech3} suppose that
one of the graphs $A$, $B$ is standard. If $\hat a_i=\hat b_j$
for some $i, j$ then $A$ can be obtained from $B$
by a sequence of turns and shifts such that after an
appropriate inner blowup of $\G$, every step is dominated by
$\G$.
\elem

\bproof By Lemma \ref{lemshift}  it suffices to distinguish the
following cases.

Case 1: $n=m=0$ and $i=2k+2$, $j=2l+2$;

Case 2: $n=0$, $i=2k+2$ and $j=1$;

Case 3: $i=j=1$.

Moreover by {\it loc.cit.}\ in the case $n\ge 1$  we may assume
that $\alpha_n\le -1$.

We note that under our assumptions the birational map
$A\dashrightarrow B$ restricts to  one between the  chains
$A':=A\ominus \{a_i\}$ and $B':=B\ominus\{b_{j}\}$ dominated by the
chain $\G':=\G\ominus\{\hat a_i=\hat b_{j}\}$.

\noindent $\bullet$ In case 1 the  chains
$A'=[[0_{2k+1}]]$ and
$B'=[[0_{2l+1}]]$ are birationally
equivalent, so by Proposition \ref{p2} $k=l$ and
$A'\dashrightarrow B'$ is equal to $\tau^s$, where $\tau$ is a
move. This amounts to an $s$-iterated shift $A\dashrightarrow B$,
and the assertion follows.

\noindent $\bullet$ In case 2 we have $\hat a_{2k+2}=\hat b_1$,
hence the birational map
$$A=((0_{2k+1},\alpha_0))\dashrightarrow
B=((0_{2l+1},\beta_0,\ldots,\beta_m))
$$
restricts to $$A'=[[0_{2k+1}]]\dashrightarrow
B'=[[0_{2l},\beta_0,\ldots,\beta_m]]\,.$$ By Remark
\ref{remzeros}(2) the standard form of the chain on the right can
have at least $2l$ and at most $2l+2$ zeros, therefore $k=l$.
Using Lemma \ref{l1}.b repeatedly we get $\hat a_{2t}=\hat
b_{2t+1}$ for $t=1,\ldots, k$.

If $\hat a_1= \hat b_2$ in $\G$ then by Lemma \ref{l1}.b, $\hat
a_{2t-1}=\hat b_{2t}$ for $t=1,\ldots, k+1$. In other words, $\hat
A\subseteq \hat B=\G$ and so $A$ is a blowdown of $B$, which
implies that $A= B$. For $A$ and $B$ are both almost standard and
$B$ does not admit a nontrivial blowdown to another almost
standard circular graph.

Suppose now that $\hat a_1\neq \hat b_2$ in $\G$. Again by Lemma
\ref{l1}.b, either $$ \hat b_{2t}\in ]\hat a_{2t-1},\hat
a_{2t}=\hat b_{2t+1}[\,\,\,\forall t=1,\ldots,k+1\,\,\,\mbox{
or}\,\,\, \hat a_{2t-1}\in ]\hat b_{2t},\hat a_{2t}=\hat
b_{2t+1}[\,\,\,\forall t=1,\ldots,k+1.$$ In the first case we
perform a shift of $A$ by blowing up the edges $[a_{2t-1},a_{2t}]$
and blowing down $a_{2t-1}$ for $t=1,\ldots k+1$. Similarly, in
the second case we perform a shift in the other direction. In both
cases $\G$ will dominate the corresponding shift $A\dashrightarrow
\tilde A$. Replacing $A$ by $\tilde A$ diminishes the distance
between $\hat a_1$ and $\hat b_2$ in $\G$ and does not affect the
vertices $\hat a_{2t}=\hat b_{2t+1}$, $t=1,\ldots, k$
(cf. Remark
\ref{remshift}(3)). After a finite sequence of shifts we can
achieve that $\hat a_1= \hat b_2$ and hence $\hat A=\hat B$. This
concludes case 2.

\noindent $\bullet$ In case 3 we may suppose that $k\le l$. Since
$\hat a_{1}=\hat b_{1}$, the map $A\dashrightarrow B$ restricts to
$$A'=[[0_{2k},\alpha_0,\ldots, \alpha_n]]
\dashrightarrow B'=[[0_{2l},\beta_0,\ldots, \beta_m]]\,.$$
Using Lemma \ref{l1}.a, $\hat a_{2s+1}=\hat b_{2s+1}$ $\forall
 s=1,\ldots, k$; this is true as well if $k=0$.
 Hence $A\dashrightarrow B$ further restricts to
 a birational map
 \be\label{as}
A''=[[\alpha_0,\ldots, \alpha_n]] \dashrightarrow
B''=[[0_{2(l-k)},\beta_0, \ldots, \beta_m]]\,.\ee Let us consider
the cases $k<l$ and $k=l$ separately.

\noindent $\diamond$ In the case $k<l$ necessarily $n\ge 1$ and
$\alpha_0=0$, so $A$ is standard, since otherwise $A''$ would
transform into a standard chain with at most one zero weight,
whereas the standard form of $B''$ has at least 2 zero weights.
Comparing the standard models of $A''$ and $B''$ yields as well
that $l=k+1$. Applying Lemma \ref{l1}.a gives $\hat a_{2k+3}=\hat
b_{2k+3}$.  By Lemma \ref{l1}.b, the order in the pair of vertices
$\hat a_2,\hat b_2$ in $\G$ is inherited by the pairs $\hat
a_{2t},\hat b_{2t}$ $\forall t=1,\ldots,k+1$. Hence as before, if
$\hat a_2\neq \hat b_2$ then $\G$ dominates a shift of $A$ which
diminishes the distance between $\hat a_2$ and $\hat b_2$. By a
finite sequence of such shifts we can transform $A$ into a graph
\be \label{asbefore} C=((\eta,0_{2k+1},
\alpha_1-\eta,\alpha_2,\ldots,\alpha_n))\,,\quad \eta\in \Z, \ee
with vertices $c_1,\ldots,c_{2k+n+2}$ such that $\hat
c_{2t-1}=\hat a_{2t-1}=\hat b_{2t-1}$ for $t=1,\ldots, k+2$ and,
moreover, $\hat c_2= \hat b_2$. Thus by Lemma \ref{l1}.b, $\hat
c_{2t}=\hat b_{2t}$ for $t=1,\ldots, k+1$.

We may further assume that the pair of dominations $\G\to C$ and
$\G\to B$ is relatively minimal i.e., $\G\ominus (\hat C\cup \hat
B)$ does not contain a $(-1)$-vertex. Since the weights
$\alpha_2,\ldots, \alpha_{n-1}$ of $C$ are $\le -2$ and $\hat
c_t=\hat b_t$ $\forall t=1,\ldots,2k+3$, no $(-1)$-vertex in $\hat
C$ can be blown down in $B$. Hence $\G=B$ and so $C$ is a blowdown
of $B$. Since the weight of a vertex can only increase under a
blowdown, this gives that $\eta,\alpha_1-\eta\ge 0$ and thus
$\alpha_1\ge 0$, which is a contradiction.

\noindent $\diamond$ Finally let us suppose that $k=l$. By
assumption, one of the graphs, say $A$, is standard, see
(\ref{cstandard}). Thus in (\ref{as})
$$
A''=[[0,\alpha_1, \ldots, \alpha_n]] \sim B''=[[\beta_0, \ldots,
\beta_m]]\,,
$$
unless $A=((0_{2k+1},w))$ and $A''=[[w]]$ with $w\le -1$.

In the latter case either $B''=\emptyset$ and interchanging $A$
and $B$ returns us to case 2, or $B''\neq \emptyset$ contracts to
$A''$ and so $B$ contracts to $A$. However since
$a_1^2=b_1^2=a_{2k+1}^2=b_{2k+1}^2=0$ no contraction in $B$ is
possible, so $A=B$, as required.

In the former case, if $n=0$ then after renumbering the
vertices of $A$ we are in case 2, which was already treated. If
$n\ge 1$ then similarly as above $\beta_0=0$ and $B$ is as well
standard, cf. Definitions (\ref{cstandard}) and \ref{tech}. By
Lemma \ref{l1}.a we get $\hat a_{2i+1}=\hat b_{2i+1}$ for
$i=0,\ldots, k+1$. Applying shifts to $A$ as in the previous case
we can again transform $A$ into a graph $C$ as in (\ref{asbefore})
which satisfies $\hat c_i=\hat b_i$ for $i=1,\ldots ,2k+3$.
Arguing as before it follows that $C$ is a blowdown of $\G=B$.

If $C=B$ then we are done. If $C\ne B$ then at least one blowdown
occurs. Since $B$ is standard, this is only possible if $n=2$,
$B=[[0_{2k+2},-1,-1]]$, and then $C$ is obtained from $B$ by
blowing down the $(-1)$-vertex $b_{2k+4}$. Hence $m=1$ and
$\alpha_1-m=0$, which gives $\alpha_1=1$. This is impossible since
$A$ is standard, so indeed $C=B$ as desired. \eproof

\bsit{\it Proof of Theorem \ref{maincirc}.}\label{circp} Given a
birational morphism $A\dasharrow B$ of standard circular graphs,
we let $\G\to A$ and $\G\to B$ be a corresponding pair of
dominations. If all weights of $A$ or of $B$ are $\le -1$ then the
assertion follows by Lemma \ref{circl3}. Otherwise, by Lemma
\ref{circl2} $A$ can be transformed by a sequence of shifts into
an almost standard circular graph $A'$ such that $\hat A'\cap \hat
B\ne \emptyset$. Applying Lemma \ref{circl5} to $A'$ and $B$ the
result follows. \qed \esit

From Theorem \ref{maincirc} we deduce the following.

\bcor\label{corcirc} Every standard circular graph
$A=((0_{2k},a_1,\ldots, a_n))$ with $k,n\ge 0$ and $a_1,\ldots,
a_n\le -2$ is unique in its birational equivalence class up to
cyclic permutation and reversion of the sequence
$(a_1,\ldots,a_n)$. Moreover the standard circular graphs
$((0_l,w))$ ($l\ge 0$, $w\le 0$) and $((0_{2k},-1,-1))$ ($k\ge
1$) are also unique in their birational equivalence classes. \ecor

\bproof If $A$ is a graph as in \ref{tech} then the sequence of
numbers $\alpha_0+\alpha_n,\alpha_1,\ldots,\alpha_{n-1}$ remains
unchanged by turns and shifts up to a cyclic permutation and
reversion. Since by Theorem  \ref{maincirc} an arbitrary
birational transformation is a composition of turns and shifts,
the result follows. \eproof

\subsection{Proof of Theorem \ref{graph.8}}

We need the following stronger form of Proposition \ref{graph.9}.

\blem \label{graph.8a} If the weighted graph $\Gamma$ dominates
two non-circular semistandard graphs $A,B$ then there is a
sequence of birational transformations
$$
\bdi[size=1.5em] &&&&\Delta_1'&&&& \Delta_2'&&\ldots &&
\Delta_n'\\
&&&\ldTo &&\rdTo&&\ldTo &&\rdTo&&\ldots && \rdTo\\
A&=&\Delta_0&&&& \Delta_1 &&&&\Delta_2&&\ldots && \Delta_n&=&B\,,
\edi\leqno (*)
$$
where all the $\Delta_i$ are semistandard and $\Delta_i'\to
\Delta_{i-1}$ and $\Delta_i'\to \Delta_{i}$ are admissible.
Moreover, after a finite sequence of inner blowups of $\Gamma$
all the $\Delta_i'$ are dominated by $\Gamma$.
\elem

\bproof By Proposition \ref{graph.9} we can find a sequence as in
$(*)$, where all the $\Delta_i'\to \Delta_{i-1}$ and $\Delta_i'\to
\Delta_{i}$ are admissible. To obtain from such a sequence another
one in which all the $\Delta_i$ are also semistandard, we proceed
as follows. By Theorem \ref{graph.7} there is a semistandard graph
$\tilde\Delta_i$ which is obtained from $\Delta_i$ by an inner
transformation. By Lemma \ref{anal0} there are sequences of inner
blowups $\tilde\Delta_i'\to \Delta_i'$ and $\tilde\Delta_{i+1}'\to
\Delta_i'$ such that $\tilde\Delta_i$ is dominated by
$\tilde\Delta_i'$ and $\tilde\Delta_{i+1}'$. Using again Lemma
\ref{anal0} after inner blowups $\G$ will dominate
$\tilde\Delta_i'$ and $\tilde\Delta_{i+1}'$. Replacing $\Delta_i$,
$\Delta_i'$, $\Delta_{i+1}'$ by $\tilde\Delta_i$,
$\tilde\Delta_i'$, $\tilde\Delta_{i+1}'$, respectively, the new
graph $\Delta_i$ will be semistandard. After a finite number of
steps we arrive at a sequence $(*)$ as required. \eproof

In general, one cannot achieve that all graphs $\Delta_i$ as in
the lemma above are standard. This can be seen by the following
simple example.

\bexa\label{nonlindom}
Consider the diagram of dominations
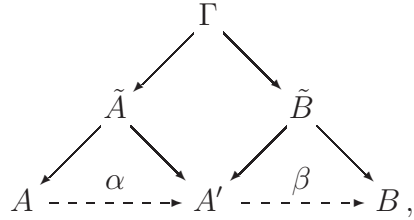
\bdi[size=1.5em]
   &&&&& \Gamma \\
   &&&&\ldTo &&\rdTo \\
   &&& \tilde A&&&&\tilde B \\
   &&\ldTo&&\rdTo &&\ldTo &&\rdTo \\
   &A &&\rDashto^{\alpha}&&A'&&\rDashto^{\beta}&&B\,,
\edi
where $\alpha,\,\beta$ are elementary transformations,
$\Gamma$ is the non-linear graph

\bigskip
$$\Gamma:\qquad\quad\quad
 \cshiftdown{-1}{b}
 \cshiftup{-1}{a}\phantom{\lin}
 \nwlin\cou{-2}{c}\swlin\co{}\lin\cou{-1}{d}\lin\cou{w}{e}
 \ldots\ldots\,,
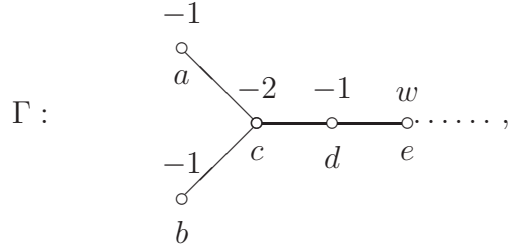
$$ \vskip 0.5in

the other graphs are the linear chains

\bigskip
$$\tilde A:\quad \quad
 \cshiftup{-1}{a}\phantom{\lin}
 \nwlin\cou{-1}{c}
 \lin\cou{-1}{d}\lin\cou{w}{e}\ldots\ldots\,,\qquad
\qquad\tilde B:\quad\quad
 \cshiftdown{-1}{b}\phantom{\lin}
 \cou{-1}{c}
 \swlin\co{}\lin\cou{-1}{d}\lin\cou{w}{e}\ldots\ldots\,,
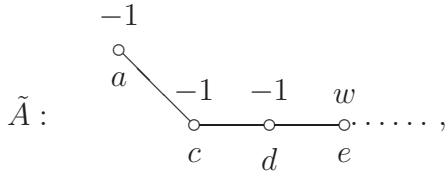
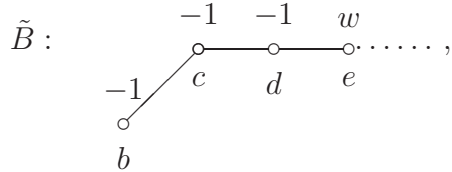
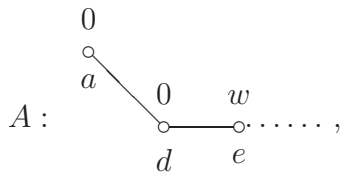
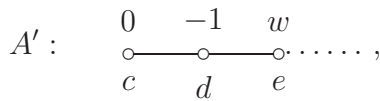
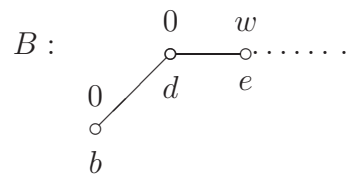
$$ \bigskip

 \bigskip
$$A:\quad
 \cshiftup{0}{a}\phantom{\lin}
 \nwlin\cou{0}{d}
  \lin\cou{w}{e}\ldots\ldots\,,\qquad
 A':\quad \quad
 \cou{0}{c}
 \lin\cou{-1}{d}\lin\cou{w}{e}\ldots\ldots\,,\qquad
B:\quad
 \cshiftdown{0}{b}\phantom{\lin}
 \cou{0}{d}\swlin\co{}\lin\cou{w}{e}\ldots\ldots\,.$$ \vskip 0.5in

 \noindent
 We notice that
 $\Gamma\to A, \,\,\Gamma\to B$ is a relatively minimal
 pair of dominations. To transform $A'$ into a standard
 chain an outer blowup at $c$ is required.
 However, it is not possible to lift this outer blowup to $\G$,
 $\tilde A$ and $\tilde B$ simultaneously in such a way that
 $\tilde A$
 and $\tilde B$ remain both linear.\eexa

 \brem\label{philo} It is worthwhile to compare this example
 with Lemma \ref{anal0}.
 Letting in this lemma $\G=B,\,\G'=\tilde B$ and $\tilde G=A$
 as in Example \ref{nonlindom}
 and
 proceeding as in the proof of the lemma
 yields a pair of dominations of $A$ and $\tilde B$ from
 a linear chain
 $\tilde \G'=\tilde B$.
 The crucial difference consists in the following.
 Attributing names
 to vertices of the graph $\G$ as in \ref{graph.8a}
 or \ref{nonlindom} gives also names to vertices of $A$ and $B$.
 This in general makes the procedure of the proof of
 \ref{anal0} impossible in the 1-st case of (iii).
 Actually keeping track of these names, if the outer blowups at
 $v_0$ gives
 two different vertices in $\tilde \G$ and in $\G'$, then
 we must create a branch point
 of $\tilde \G'$. Fortunately, under the setup as in Lemma
 \ref{anal0}, no names are prescribed in advance. \erem

\bproof[Proof of Theorem \ref{graph.8}.] Let $A\dashrightarrow B$
be a birational map between two semistandard graphs as in the
theorem. By Lemma \ref{graph.8a} it suffices to show that for any
two admissible dominations $\sigma_i:\Gamma\to \G_i$, where $\G_i$
are semistandard graphs ($i=1,2$), $\gamma=\sigma_2\circ
\sigma_1^{-1}$ can be decomposed into a sequence of elementary
transformations dominated by $\Gamma$ after, possibly, some inner
blowups of $\Gamma$.

Since $\sigma_i$ is admissible, $B(\G_i)\hat{}=B(\G)$, $i=1,2$,
and so for every segment $\Sigma$ of $\G$, $\sigma_i$ induces an
admissible domination $\sigma_i\vert\Sigma:\Sigma\to \Sigma_i$
onto the contraction $\Sigma_i$ of $\Sigma$ in $\G_i$. It suffices
to show that $\gamma\vert\Sigma_1:\Sigma_1\dasharrow \Sigma_2$ can
be decomposed into a sequence of elementary transformations
dominated by some inner blowup of $\Sigma$. The latter follows
from Propositions \ref{c1} and \ref{circl3}\footnote{See also
Corollary \ref{cornew}.} for semistandard linear segments
$\Sigma_1,\, \Sigma_2$ and for standard circular graphs,
respectively. \eproof

We deduce the following analog of Corollary \ref{corcirc}.

\bcor\label{c2} The standard form of a linear chain is unique in
its birational class up to reversion. \ecor

\bproof If $\G\to A$, $\G\to B$ is a pair of dominations of
standard linear chains $A$, $B$, then we can factor the
corresponding birational transformation $A\dasharrow B$ as in
Lemma \ref{graph.8a}. Using Remark \ref{sstos}(1), each  of the
semistandard graphs $\Delta_i$ in {\it loc.cit.}\ can be
transformed into a standard linear chain, say $A_i$, $1\le i\le
n-1$, by a sequence of elementary transformations. It suffices to
verify that $A_i$ is equal to $A_{i+1}$ or to its reversion for
$i=0,\ldots, n$. Using Lemma \ref{anal0} the birational map
between $A_i$ and $A_{i+1}$ can be dominated by some linear chain,
say, $L_i$\footnote{In general, neither $A_i$ nor $L_i$ are
dominated by $\G$.}. In other words, it suffices to treat the
case where $A$ and $B$ are dominated by a linear graph.

In this situation, by virtue of Proposition \ref{p1}, $A=B$ or
$A=B^*$  as soon as at least one of the standard chains $A$ and
$B$ with $A\sim B$ is different from $[[0_{2k+1}]]$ for all $k$.
Otherwise $A=[[0_{2k+1}]]$ and $B=[[0_{2l+1}]]$, and Proposition
\ref{p2} yields that $k=l$. Alternatively, $k$ and $l$ are equal
to the number of positive eigenvalues of the associate bilinear
form $I(A)$, $I(B)$, respectively, and this number is a birational
invariant, see \ref{adjma} in the Appendix below.  \eproof

Finally, Corollaries \ref{c2} and \ref{corcirc}
yield the following result.

\bcor\label{unique}
Every non-circular standard weighted graph is unique in its birational
equivalence class up to reversions of its linear segments. Similarly,
a circular standard graph is unique in its birational equivalence class up
to a cyclic permutation of its nonzero weights and reversion.
\ecor

\subsection{The geometric meaning}
In the geometric setup, $\G$ is the dual
graph of a reduced divisor $D$ with normal crossings (an NC-divisor, for short)
on the regular part $X_{\rm reg}$ of a normal complete algebraic surface $X$.
If $V=X\setminus D$ then any two such completions
$V\hookrightarrow X_1, \,V\hookrightarrow X_2$ can be dominated by a third one
$V\hookrightarrow X$:
\be\label{diagr1}
\bdi[size=1.5em]
&&& X &&&\\
&&\ldTo^{} &&\rdTo^{}&& \\
& X_1& & \rDashto & &X_2&
\edi
\ee
inducing the identity on $V$. Letting $\Gamma_i,\,i=1,2$, and $\Gamma$
be the dual graphs of the corresponding boundary divisors $D_1,D_2$ and $D$,
respectively, we get a diagram of dominations
\be\label{diagr2}
\bdi[size=1.5em]
&&& \Gamma &&&\\
&&\ldTo^{} &&\rdTo^{}&& \\
& \Gamma_1& & \rDashto & &\Gamma_2&. \edi \ee

\bprop\label{geo} Suppose we are given a diagram (\ref{diagr2})
and two SN-completions $(X_1,D_1)$ and $(X,D)$ of the same normal
algebraic surface $V$ such that $(X,D)$ dominates $(X_1,D_1)$ and
induces the given domination $\Gamma\to\Gamma_1$ for the dual
graphs of $D$ and $D_1$, respectively. Assume moreover that $\G\to
\G_2$ contracts only vertices corresponding to rational curves of
$D$. Then there exists a unique SN-completion $(X_2,D_2)$ of $V$
dominated by $(X,D)$ which fits the diagram (\ref{diagr1}) and
induces the given domination $\Gamma\to\Gamma_2$ for the dual
graphs of $D$ and $D_2$, respectively. \eprop

\bproof To get $X_2$ from $X$, it is enough to contract
all the irreducible components in $D$ which correspond to the vertices in
$\Gamma\ominus\hat \Gamma_2$.
\eproof

\bdefi\label{geost} We say that an  SN-completion $(X,D)$ of $V$ is {\it
standard} if the dual graph $\Gamma$ of $D$ is. \edefi

\bcor\label{geodg} Every normal algebraic surface $V$ admits a standard
completion
$(X,D)$.
Moreover, any other such completion can be obtained from $(X,D)$ by
a sequence of elementary transformations. \ecor

\section{Appendix}
\subsection{The adjacency matrix and the discriminant of a weighted graph}
\label{adjma} For the sake of simplicity, we restrict in this
section to weighted graphs with integral weights without loops and
multiple edges. To such a graph $\Gamma$ one usually associates
its adjacency matrix $I(\Gamma)=(v_i.v_j)$, where $v_i^2=w_i$ is
the weight of the vertex $v_i$, whereas for $i\ne j$, $v_i.v_j=1$
if $v_i$ and $v_j$ are joined by an edge, and $v_i.v_j=0$
otherwise. It is well known \cite[Prop. 1.1]{Ne}, \cite[Prop.
1.14]{Ru} that a blowup of $\Gamma$ just adds a negative
eigenvalue to this matrix, see (\ref{ma}) below. In particular,
the number of positive (non-negative) eigenvalues is a birational
invariant of $\Gamma$, and so is its {\em discriminant} $\delta
(\G)=\det (-I(\G))$. If $\Gamma'$ is a subgraph of $\Gamma$ then
$I(\Gamma')$, being a symmetric submatrix of $I(\Gamma)$, has at
most the same number of positive eigenvalues as $I(\Gamma)$.

\bsit\label{deter} For a weighted graph $\G$ and a vertex $v$ in
$G$ of weight $a=v^2$, we let
$$
\delta_v (\G)= \delta(\G\ominus \{v\})=\prod_j \delta(\G_j)\,,
$$
where $\G_j$ runs through the set of branches of $\G$ at $v$.  If
$\G_j$ is joined to $v$ by a unique edge $[v,v_j]$ then the
following
  holds (see e.g., \cite{DrGo, Ne, OrZa}):
\be\label{grdi} -\delta(\G)=a\cdot \delta_v(\G)+\sum_{j}
\delta_{v_j} (\G_j)\cdot\prod_{i\neq j} \delta(\G_i)\,. \ee \esit

\bsit\label{contt} Further, if $\G$ is contractible i.e.,
dominates the one-vertex graph with weight $-1$ then clearly $\G$
is a tree and $I(\G)$ is negative definite of discriminant $\delta
(\G)=1$. Vice versa, we have the following lemma. \esit

\blem\label{graph.3}
\begin{enumerate} [(a)]\item  A weighted tree $\G$ with
negative definite adjacency matrix $I(\G)$  and discriminant
$\delta (\G)=1$ is contractible.
\item Consequently, a weighted graph $\G$ with
discriminant $1$ that can be transformed into a one-vertex graph
is contractible.
\end{enumerate}
\elem

\bproof For (a) we refer the reader to \cite{Mu, Hi2} or
\cite[Prop. 1.20]{Ru}. (b) is an immediate consequence of (a) due
to the fact that the discriminant does not change under
birational transformations. \eproof

The following lemma is well known, see e.g., \cite[3.8]{Fu} or
\cite[3.3.1]{Mi}.\esit

\blem\label{cfr} Let $\frac{m}{e}$ be the continued fraction
$$\frac{m}{e}=[k_1,\ldots,k_n]=k_1-\frac{1}{k_2-\frac{1}{\ddots
-\frac{1}{k_n}}}\qquad\mbox{with}\quad 0\le e
<m,\,\,\,\gcd(e,m)=1\,,$$ where $k_i\in\N_{\ge 2}\,\,\forall i\ge
2\,,$ and let $\G$, $\G'$ be the linear chains $[[-k_1,\ldots,-k_n]]$,
$[[-k_2,\ldots,-k_n]]$, respectively. Then $m=\delta(\G)$ and
$e=\delta(\G')$. \elem

Thus the pair $\left(\delta(\G),\delta(\G')\right)$
uniquely determines the
linear chain $\G$.

\bsit\label{proce} We let $L=[[w_1,\ldots,w_n]]$ ($n\ge 1$) be a linear
branch of a weighted graph $\G$ at a vertex $v_0\in\Gamma^{(0)}$, where
$w_1$ is the weight of the neighbor $v_1$ of $v_0$ in $L$ and
$w_n$ is the weight of the end vertex $v_n$ of $L$. We let $\G_0$ be the
graph obtained from $\G$ by deleting the branch $L$ and
changing the weight $w_0=v_0^2$ to
$$w'_0=
-[-w_0,-w_1,\ldots,-w_n]=w_0+[-w_1,\ldots, -w_n]^{-1}=
w_0+\frac{1}{w_1+\frac{1}{\ddots
+\frac{1}{w_n}}}\,. $$
We denote
by $ \langle a\rangle$ the usual Euclidean quadratic form on
$\R^1$ multiplied by $a$. \esit

\blem\label{cfr1} In the notation as above suppose that
$w_1,\ldots ,w_{n-1}\le -2$ and $w_n\le -1$. Then
\be\label{byind} I(\G)\sim
I(\G_0) \oplus \langle w_1'\rangle\oplus\ldots \oplus \langle
w_n'\rangle\,\qquad\mbox{and}\qquad \delta (\G)=\delta (\G_0)\cdot
\prod_{i=1}^n (-w_i')\,, \ee where
\be\label{byind1}
w_{i}'=-[-w_i,\ldots,-w_n]\le -1\,\qquad\forall i=1,\ldots,n\,.
\ee Moreover, $w_{i}'<-1$ $\forall i=1,\ldots,n$ if also $w_n\le-2$.
Consequently, $\G$ and $\G_0$ have the same number of positive
(non-negative) eigenvalues. \elem

\bproof Essentially, our proof repeats an argument in \cite[p.
20]{Mu}. We proceed by induction on $n$. In the case $n=0$ the linear
branch $L$ is empty and the assertion is immediate. Assume
now that $n\ge 1$.  In the linear space $V(\G)$ spanned over $\R$ by
the vertices of $\G$, we consider the symmetric bilinear form defined by
the adjacency matrix $I(\G)$. The vector
$v_{n-1}'=v_{n-1}-w_n^{-1}v_n$ is orthogonal to $v_n$ and
satisfies
$$
v_{n-1}^{\prime2}=v_{n-1}^2- \frac{1}{w_n}, \qquad
v_{n-1}'v=v_{n-1}v \mbox{ for all vertices of $\Gamma$ different
from }v_{n-1}, v_n.
$$
Hence deleting  the vertex $v_n$ from $\Gamma$ and changing the
weight of $v_{n-1}$ to $v_{n-1}^{\prime 2}=w_{n-1}'$ with
$w_{n-1}'=w_{n-1}-1/w_n$, we obtain a new graph $\Gamma_{n-1}$
which satisfies $I(\G)\sim I({\G_{n-1}})\oplus \langle
w_n\rangle$. Since $w_{n-1}'\le -1$, this completes the induction
step. \eproof

\bexas\label{exex} 1. If $L=[[-2,\ldots,-2,-1]]$ (of length $n$) then
$w_0'=w_0+1$ and $I(\G)=I({\G_0})\oplus \langle
-1\rangle\oplus\ldots\oplus \langle -1\rangle$ ($n$ times).

2. If $L=[[-1,-2,\ldots,-2]]$ (of length $n$) then $w_0'=w_0+n$   and
$$I(\G)\sim I({\G_0}) \oplus \left\langle  -\frac{1}{n}
\right\rangle \oplus\left\langle -\frac{n}{n-1}\right\rangle
\oplus\ldots\oplus \left\langle -\frac{r}{r-1}\right\rangle
\oplus\ldots\oplus \left\langle -\frac{3}{2}\right\rangle \oplus
\left\langle -2\right\rangle \,.$$ In both cases $\delta
(\G)=\delta (\G_0)$. \eexas

\brems\label{2prop} (1) In the geometric setting, $\G$ is the dual
graph of a reduced divisor $D$ on a normal surface $X$ and $\G_0$
is the dual graph of the image of $D$ under contraction of the
part $D'$ of $D$ with dual graph $\G'$. The orthogonal basis of
the reduction as in the proof above appears geometrically as
follows. For a given $i$ with $0\le i< n$ we consider the
contraction $\sigma_i: X\to X_i$ of all the components $C_j$ of
$D'$ with $j> i$. The total transforms $\sigma_i^* (\sigma_i
(C_i))$, $i=1,\ldots,n$, on the original surface $X$ are then
mutually orthogonal due to the projection formula. Moreover all of
them are orthogonal to the total transform $\sigma_0^* (\sigma_0
(\G))$ which has dual graph $\G_0$.

(2) More generally, we let $\Gamma\to\Gamma_1$ be a domination of
weighted graphs consisting of a sequence of blowdowns
$$\Gamma=\Gamma_{n+1}\to\Gamma_{n}\to\ldots\to\Gamma_1\,.$$
Letting also $v_1,\ldots,v_k$ be the vertices of $\Gamma_1$ and
$u_i$ be the blowdown vertices in
$\Gamma_{i+1}\ominus\widehat{\Gamma_i}$, $i=1,\ldots,n$, we
consider as before the vector space $ V(\Gamma)$ endowed with the
symmetric bilinear form $I(\Gamma)$. It is easily seen that the
subspace $V_1={\rm span}\,(v_1^*,\ldots,v_k^*)$ has the orthogonal
complement $V_1^\bot={\rm span}\,(u_1^*,\ldots,u_n^*)$, where
$u_1^*,\ldots,u_n^*$ are mutually orthogonal eigenvectors of
$I(\Gamma)$ corresponding to negative eigenvalues. Moreover the
restriction of $I(\Gamma)$ to $V_1$ is equivalent to the bilinear
form $I(\Gamma_1)$ on $V(\Gamma_1)$. Therefore \be\label{ma}
I(\Gamma)\sim I(\Gamma_1)\oplus\langle-1\rangle^n\,.\ee \erems

\bprop\label{3prop} For a linear chain $L=[[w_0,\ldots,w_n]]$
with integral weights the following hold.
\begin{enumerate}[(a)]\item If $w_i\le -2$
$\forall i=0,\ldots,n$ then $I(L)$ is negative definite of
discriminant $\delta (L)>1$. \item   If $w_{k}=-1$ for some
$k\in\{0,\ldots,n\}$ and $w_i\le-2$ $\forall i\neq k$ then $L$
is contractible if and only if $\delta (L)=1$ or, equivalently,
\be\label{fla} \frac{e_1}{m_1} + \frac{e_2}{m_2}
=1-\frac{1}{m_1m_2}\,,\ee where
$$\frac{e_1}{m_1}=[-w_{k-1},\ldots,-w_0]^{-1}\quad\mbox{and}\quad
\frac{e_2}{m_2}=[-w_{k+1},\ldots,-w_n]^{-1}\,$$ with $m_i>0$ and
$\gcd (e_i,m_i)=1$, $i=1,2$.
\item
If $w_0=0$, $w_1\in\Z$ and $w_i\le -2$ $\forall i=2,\ldots,n$, or $n=2$,
$w_0=0$, $w_1\in\Z$ and $w_2\le 0$ then $I(L)$ has exactly one positive
eigenvalue, and $\delta (L)<-1$ if $n\ge 3$, $\delta (L)\le 0$
if $n=2$. \end{enumerate}\eprop

\bproof (a) follows from Lemma \ref{cfr1} above applied to
$\G_0=v_0$.

Similarly, to show (b) we apply the reduction as in Lemma
\ref{cfr1} to the two branches $L_1=[[w_0,\ldots,w_{k-1}]]$ and
$L_2=[[w_{k+1},\ldots,w_n]]$ at the only $(-1)$-vertex $v_k$ of
$\G$. By virtue of Lemma \ref{cfr1}, $\G$ and $\G_0=[[w_k']]$ have
the same number of positive (non-negative) eigenvalues. Hence
$\G$ is negative definite if and only if $w_k'<0$. Using
(\ref{proce})
$$
w_k'=-1+[-w_{k-1},\ldots,-w_0]^{-1}+[-w_{k+1},\ldots,-w_n]^{-1}=
-1+e_1/m_1+e_2/m_2\,,
$$
where by Lemma \ref{cfr}, $m_i=\delta (\G_i)>0$ and $e_i=\delta
(\G_i\ominus v_{k\pm 1})\ge 0$. On the other hand, according to
(\ref{grdi}) and (\ref{byind}),
$$\delta (\G)=(-w_k') \delta (\G_1)\delta
(\G_2)=m_1m_2-m_1e_2-m_2e_1=1  \\
\,\,\,\Longleftrightarrow\,\,\,
-w_k'=1-\frac{e_1}{m_1}-\frac{e_2}{m_2} =\frac{1}{m_1m_2}\,.$$
Now (b) is a consequence of Lemma \ref{graph.3}.a.

In case (c), letting $\langle
a,b\rangle=I([[a,b]])=\begin{pmatrix}
   a & 1 \\
   1 & b
\end{pmatrix}$ and applying Lemma \ref{cfr1} to
$\G_0=[[w_0,w_1']]=[[0,w_1']]$ we obtain
$$I({\G})\sim \langle 0,w_1'\rangle\oplus \langle
w_2'\rangle\oplus\ldots \oplus \langle w_n'\rangle\sim \langle
0,0\rangle\oplus \langle w_2'\rangle\oplus\ldots \oplus \langle
w_n'\rangle\,,$$ where as before $w_k'=-[-w_{k},\ldots,-w_n]<-1$
$\forall k=2,\ldots,n$ if $n\ge 3$ and $w_2'=w_2\le 0$ if $n=2$.
Hence $\delta (\G)=-(-w'_2)\cdot\ldots\cdot(-w_n')<-1$ if $n\ge
3$ and $\delta (\G)\le 0$ if $n=2$. Thus $I(\G)$ has exactly one
positive eigenvalue, as stated. \eproof

\bexas\label{tob} According to (b),
the equality $\frac{1}{2}+\frac{1}{3}=1-\frac{1}{6}$
leads to the contractible linear
chain $L=[[-2,-1,-3]]$. Similarly,
$\frac{2}{3}+\frac{1}{4}=1-\frac{1}{12}$ yields
$L=[[-2,-2,-1,-4]]$.
\eexas

\subsection{Spectra of linear and circular standard graphs}

The spectrum and the inertia indices $i_\pm, i_0$ of a weighted
graph $\Gamma$ are as usual the spectrum, respectively, the
inertia indices of the associate symmetric bilinear form
$I(\Gamma)$. In the same way as it was done in the proofs of Lemma
\ref{cfr1} and Proposition \ref{3prop}.c, we can find the
spectra of the standard graphs. Let us start with the standard
linear chains.

\bprop\label{numeig} \begin{enumerate}\item[(a)] The form
$I([[0_{m}]])$ has $\lfloor\frac{m}{2}\rfloor$ negative and
$\lfloor\frac{m}{2}\rfloor$ positive eigenvalues; for an odd $m$
it has in addition a zero eigenvalue. Thus the inertia indices of
$\Gamma=[[0_{m}]]$ are
$$(i_+,i_-,i_0)=\begin{cases} (k,k,0), & m=2k\\
(k,k,1), & m=2k+1\,\end{cases}$$ and the discriminant is $\delta
([[0_{2k}]])=(-1)^k$ and $\delta ([[0_{2k+1}]])=0$ $\forall k\ge
0$. \item[(b)] If $L=[[0_{2k},w_1,\ldots,w_n]]$ and
$L'=[[w_1,\ldots,w_n]]$ then
$$I(L)\sim I(L')\oplus I([[0_{2k}]])\sim I(L')
\oplus\bigoplus_{i=1}^k\langle 0,0\rangle=I(L') \oplus\langle
0,0\rangle^k\,.$$ Consequently, if $w_i\le -2$ $\forall
i=1,\ldots,n$ then the intersection form
$$I(L)\sim \langle w_1'\rangle\oplus \langle w_2'\rangle
\oplus\ldots \oplus \langle w_n' \rangle\oplus I([[0_{2k}]])\,,$$
where $w_i'$ are defined as in (\ref{byind1}), has $n+k$ negative
and $k$ positive eigenvalues, and $\delta (L)=(-1)^k\delta (L')$.
\footnote{For the latter equality, see also
\cite[L.4.12(1)]{Dai2}.}
\end{enumerate}
\eprop

Similarly, for circular graphs the following hold.

\blem\label{cirform}\begin{enumerate}\item[(a)]
The eigenvalues $\lambda$
of a circular graph $\Gamma=((0_m))$ are
$$\left\{2\cos \frac{2\pi
l}{m}\mid l=0,\ldots,m-1\right\}\,.$$
All these eigenvalues have
multiplicity $2$ except for $\lambda=\pm 2$ if $m$ is even and
$\lambda=2$ if $m$ is odd, the latter ones being simple.
\item[(b)] For $\Gamma=((0_{4k},w_1,\ldots,w_n))$ with $n\ge 1$
we let $\Gamma'=((w_1,\ldots,w_n))$ if $n\ge 2$ and
$\Gamma'=[[w_1+2]]$ if $n=1$. Then $$I(\Gamma)\sim
I(\Gamma')\oplus\bigoplus_{i=1}^k
I\left([[0_4]]\right)=I(\Gamma')\oplus\langle0,0,0,0\rangle^k\,.$$
Consequently the inertia indices
of $\Gamma$ and $\Gamma'$ are related via
$$(i_+,i_-,i_0)=(i'_++2k,\,i'_-+2k,\,i'_0)\,.$$
\end{enumerate}
\elem

\bproof To show (a) we notice that  $I(\Gamma)$ is the matrix of
the linear map $\tau+\tau^{-1}$, where $\tau : v_i\longmapsto
v_{i+1},$ $i=1,\ldots,m$ is the cyclic shift acting on
$V(\Gamma)$. Now the statement follows from the Spectral Mapping
Theorem. Indeed $\tau^m=$id and so the complex spectrum of $\tau$
consists of the $m$-th roots of unity. Hence
$${\rm spec}\,(\tau+\tau^{-1})=\{\lambda+\lambda^{-1}\mid
\lambda\in\C,\,\lambda^m=1\}\,.$$

To show (b), in the vector space $V(\Gamma)=\bigoplus_{i=1}^m \R v_i$,
where $m=4k+n$, we perform
the following base change:
$$(v_1,\ldots,v_m)\longmapsto
(v_1,\ldots,v_4,v_5+v_1-v_3,v_6,\ldots,v_{m-1},v_m-v_2+v_4)
\qquad\mbox{if}\quad m\ge 6\,,$$ or
$$(v_1,\ldots,v_5)\longmapsto
(v_1,\ldots,v_4,v_5+v_1-v_2-v_3+v_4)\qquad\mbox{if}\quad
m=5\,\,.$$ It is easily seen that for $k\ge 1$, in the new basis
the intersection form $I(\Gamma)$ coincides with those of the
disjoint union of the graphs $[[0_4]]$ and
$((0_{4(k-1)},w_1,\ldots,w_n))$ if $m\ge 6$, respectively
$[[0_4]]$ and $[[w_1+2]]$ if $m=5$. This allows to single out an
orthogonal direct summand $I([[0_4]])$ of $I(\Gamma)$ and so
provides a reduction from $k$ to $k-1$. After $k$ steps we obtain
the desired decomposition. Now the second assertion in (b) follows
by virtue of (a). Indeed, according to Proposition \ref{numeig}.a
the inertia indices of the bilinear form $\langle 0_4\rangle^k$
are $(i_+,i_-,i_0)=(2k,2k,0)$. \eproof

\bprop\label{cirform1} The inertia indices $(i_+,i_-,i_0)$ of the
standard circular graphs\footnote{See (\ref{cstandard}).}
$$
((0_l, w)),\qquad ((0_{2k}, -1,-1))\quad\mbox{and}\quad
((0_{2k}, w_1,\ldots, w_n))\,,
$$
where $w\le 0$, $k,l\ge 0$, $n\ge 2$ and $w_i\le -2\,\,\,\forall
i$, are given in the following tables:

$$
\noindent\begin{tabular}{|c|c|c|c|c|}
  \hline
  $\Gamma$        &
  $((0_{4l}))$    &
  $((0_{4l+1}))$  &
  $((0_{4l+2}))$  &
  $((0_{4l+3}))$
                  \\
 \hline
 $(i_+,i_-,i_0)$  &
 $(2l-1,2l-1,2)$  &
 $(2l+1,2l,0)$    &
 $(2l+1,2l+1,0)$  &
 $(2l+1,2l+2,0)$
                  \\
 \hline
\end{tabular}
$$

\bigskip

$$
\noindent\begin{tabular}{|c|c|c|c|c|}
  \hline
  $\Gamma$                    &
  $((0_{2k},w))$              &
  $((0_{4l},-1))$             &
  $((0_{4l+1},w))$            &
  $((0_{4l+3},w))$            \\
  \hline $(i_+,i_-,i_0)$      &
  $(k,k+1,0)$                 &
  $(2l+1,2l,0)$               &
  $(2l+1,2l+1,0)$             &
  $(2l+1,2l+2,1)$            \\
  \hline
  except for                  &
  $((0_{4l},w)),\,-2\le w\le 0$&
  -                           &
  -                           &
  $w=0$                       \\
  \hline
\end{tabular}
$$

\bigskip

$$
\noindent\begin{tabular}{|c|c|c|c|c|}
 \hline
 $\Gamma$                    &
 $((0_{4l},-1,-1))$          &
 $((0_{4l+2},-1,-1))$        &
 $((0_{2k},w_1,\ldots,w_n))$ &
 $((0_{4l},(-2)_n))$        \\
 \hline
 $(i_+,i_-,i_0)$             &
 $(2l+1,2l+1,0)$             &
 $(2l+1,2l+3,0)$             &
 $(k,k+n,0)$               &
 $(2l,2l+n-1,1)$             \\
 \hline
 except for                  &
  -                          &
  -                          &
 $((0_{4l},(-2)_n))$         &
  -                          \\
  \hline
\end{tabular}
$$

\bigskip

 \eprop

\bproof For the circular graphs $\Gamma=((0_m))$ the result
follows by virtue of Lemma \ref{cirform}.a or, alternatively,
\ref{cirform}.b by an easy computation. In the other cases,
according to Lemma \ref{cirform}.b it is enough to consider
graphs with at most 3 zeros.

\smallskip

\noindent $\bullet$ If $\Gamma=((w))$ then $I(\Gamma)=\langle
w+2\rangle$ and so
$$(i_+,i_-,i_0)=\begin{cases} (0,1,0), & w\le
-3\\(0,0,1), & w=-2\\(1,0,0), & w=0\,\,\mbox{or}\,\,w=-1\,.
\end{cases}
$$

\smallskip

\noindent $\bullet$ For the graphs $\G=((-1,-1))$ and
$\Gamma=((0,w))$, where $w\le 0$, we have $\det I(\G)<0$ and so
$(i_+,i_-,i_0)=(1,1,0)$.

\smallskip

\noindent $\bullet$ For graphs $\Gamma_3=((0_2,w))$ and
$\Gamma_4=((0_2,w_1,w_2))$, the respective base changes
$$(v_1,v_2,v_3)\longmapsto (v_1,v_2,v_3-v_1-v_2)\,,$$
and
$$(v_1,v_2,v_3, v_4)\longmapsto
(v_1,v_2,v_3-v_1,v_4- v_2)$$ yield decompositions
$$I(\Gamma_3)\sim \langle 0,0\rangle\oplus\langle w-2\rangle\quad\mbox{and}\quad
I(\Gamma_4)\sim \langle 0,0\rangle\oplus\langle w_1\rangle
\oplus\langle w_2\rangle\,.
$$
Hence their inertia indices are
$$
(i_+,i_-,i_0)=(1,2,0)\quad\mbox{and}\quad (i_+,i_-,i_0)=
\begin{cases} (1,3,0), & w_1,w_2\le -1\\
(1,2,1), & w_1=0,w_2\le -1\\(1,1,2), & w_1=w_2=0\,,
\end{cases}
$$
respectively.

\noindent $\bullet$ If $\Gamma=((w_1,\ldots,w_n))$ is a circular
graph with $n\ge 2$ and $w_i\le -2$ $\forall i$ then for a vector
$\vec v\in V(\Gamma)$ with coordinates $(x_1,\ldots,x_n)$,
$$I(\Gamma)\vec v .\vec v= \sum_{i=1}^n w_ix_i^2+2\sum_{i=1}^n
x_ix_{i+1} =\sum_{i=1}^n (2+w_i)x_i^2 -\sum_{i=1}^n
(x_i-x_{i+1})^2\,,$$ where $x_{n+1}:=x_1$. Thus $I(\Gamma)$ is
negative definite and so $(i_+,i_-,i_0)=(0,n,0)$ provided that
$\Gamma\ne ((\,(-2)_n\,))$. Clearly $(i_+,i_-,i_0)=(0,n-1,1)$ for
the Cartan matrix of $\Gamma=((\,(-2)_n\,))$.

\smallskip

\noindent $\bullet$  Finally in the case where $\G=((0_2,
w_1,\ldots, w_n))$, $n\ge 3$ and $w_i\le -2\,\,\forall i$, we
perform the base change
$$
(v_1,\ldots, v_{n+2})\longmapsto (v',v'', v_1'\ldots, v_n'):=
(v_1,v_2,v_3-v_1,v_4,\ldots,v_{n+1},v_{n+2}-v_2)\,.
$$
Since $v'$ and $v''$ are perpendicular to $v_i'$ for all $i$, this
results in a decomposition
$$
I(\G)\sim\langle 0,0\rangle\oplus I_n'\,.
$$
To compute $I_n'$, we note first that for $i<j$
we have $v_i'v_j'=v_{i+2}v_{j+2}$ unless $i=1$ and $j=n$,
where $v_1'v_n'=-1$. Thus similarly as before
for a vector $\vec{v}=\sum_{i=1}^n x_iv_i'$
$$I'_n\vec v .\vec v
=\sum_{i=1}^n (2+w_i)x_i^2 -\sum_{i=1}^n
(x_i-x_{i+1})^2-4x_1x_n\,.$$ Since
$(x_1-x_n)^2+4x_1x_n=(x_1+x_n)^2$ the form $I'_n$ is negative
definite and $(i_+,i_-,i_0)=(1,n+1,0)$.

\smallskip

\noindent These computations combined with Lemma \ref{cirform}.b
cover all possible cases in our tables. Now the proof is
completed. \eproof

\subsection{Zariski's Lemma}
We recall the following classical fact, see e.g. \cite[Lemma
2.11.1]{Mi}. For the sake of completeness we provide a short
argument.

\blem\label{zar} ({\rm Zariski's Lemma}) If $\Gamma\to A=[[0]]$ is
a domination and $\Gamma'\subseteq \Gamma$ is a proper subgraph
then the intersection form $I(\Gamma')$ is negative definite.
\elem

\bproof Letting $a$ be the only vertex of $A$, the null space of
the quadratic form $I(\G)\sim\langle 0\rangle\oplus\langle
-1\rangle^{l-1}$ on $V(\G)$ is $\R a^*$, where $a^*$ is the total
transform of $a$ and $l=|\G|$, see (\ref{ma}). Since $a^*$ is
supported by the whole $\G$ we have $V(\G')\cap \R a^*=\{0\}$.
Hence the quadratic form $I(\G')=I(\G)\vert V(\G')$ is negative
definite.\eproof

\brem\label{ce} It is not true in general that the number of
non-negative eigenvalues of the intersection form of a graph must
strictly diminish when passing to a proper subgraph. Consider for
instance the following dominations:

\vskip 0.3in

$$
 \Gamma:\qquad
 \co{0}\lin\cou{\phantom{-m}-1}{}\nlin\cshiftup{-1}{}
\lin\co{0}\qquad\longrightarrow\qquad A:\qquad
\co{0}\lin\co{0}\lin\co{0}\,
$$

\noindent and

$$
 \Gamma':\qquad
 \co{0}\lin\co{-1}\lin\co{0}\qquad\longrightarrow\qquad
A':\qquad
 \co{1}\lin\co{1}\quad\,.
$$
\smallskip

\noindent For all 4 graphs $A,\,\G,\,A'$ and $\G'$ we have
$(i_+,i_0)=(1,1)$, whereas $\G'$ is a proper subgraph of $\G$.
\erem

However the following partial extension of Zariski's Lemma holds.

\blem\label{zar1} Let $\Gamma\to A=[[0_{2k+1}]]$ be a domination,
where $\G$ is a linear chain. If a proper subgraph $\Gamma'$ of $
\Gamma$ is either connected or satisfies $\vert\Gamma\ominus
\Gamma'\vert>k$ then the intersection form $I(\Gamma')$ has at
most $k$ non-negative eigenvalues.\elem

\bproof For $A=[[0_{2k+1}]]$ with vertices $a_0,\ldots,a_{2k}$ the
null space of the intersection form $I(A)$ is generated by the
vector
$$\vec v=a_0-a_2+a_4-\ldots\pm a_{2k}\,.$$
Since $I(\Gamma)\sim I(A)\oplus \langle -1\rangle^l$ with
$l=|\Gamma|-2k-1$, see (\ref{ma}), $I(\Gamma)$ has inertia indices
$(i_+,i_-,i_0)=(k,k+l,1)$ and the null space $\R{\vec v}^*$.

Suppose on the contrary that $I(\Gamma')=I(\G)\vert V(\G')$ has
more than $k$ non-negative eigenvalues. Then ${\vec v}^*\in
V(\G')$, so all vertices of $\Gamma$ besides possibly $\hat
a_{2i+1}$, $i=0,\ldots,k-1$ must be in $\Gamma'$. Thus
$|\Gamma\ominus\Gamma'| \le k$ and so, by our assumptions,
$\Gamma'$ is connected. By the same argument as above it must
contain the vertices $\hat a_0,\hat a_{2k}$, hence also $\hat
a_1,\ldots,\hat a_{2k-1}$ since $\G$ is linear. Therefore
$\G'=\G$, contradicting the assumption that
$\Gamma'\subseteq\Gamma$ is a proper subgraph. \eproof

\end{document}